\documentclass[a4paper,11pt]{article}

\usepackage[margin=3cm]{geometry}
\usepackage[english]{babel}
\usepackage[dvips]{epsfig}
\usepackage{amsfonts}
\usepackage{subfigure}
\usepackage{psfrag}
\usepackage{theorem}

{\theoremstyle{plain}
\newtheorem{theorem}{Theorem}
\newtheorem{proposition}[theorem]{Proposition}
\newtheorem{corollary}[theorem]{Corollary}
\newtheorem{lemma}[theorem]{Lemma}

{\theorembodyfont {\rmfamily}
\newtheorem{example}[theorem]{Example}

\newtheorem{remark}[theorem]{Remark}

}}

\begin {document}

\title{A characterization of infinite smooth Lyndon words\thanks{with the support of {\it Fonds de recherche sur la nature et les technologies} (Qu\'ebec, Canada)}}
\author{Genevi\`eve Paquin\\
Laboratoire de math\'ematiques,  Universit\'e de Savoie, CNRS UMR 5127\\
73376 Le Bourget du Lac, France\\
{\tt Genevieve.Paquin@univ-savoie.fr}}

\maketitle

\newcommand{\QED}{\rule{1ex}{1ex} \par\medskip}

\newcommand{\nl}{\par\medskip\noindent}

\newcommand{\Proof }{{\nl\it Proof.\ }}

\newcommand{\Pref}{\hbox{\rm Pref}}
\newcommand{\Suff}{\hbox{\rm Suff}}

\newcommand{\N}{\mathbb N}

\newcommand{\A}{\mathcal{A}}

\begin{abstract} In a recent paper, Brlek, Jamet and Paquin  showed that some extremal infinite smooth words are also  infinite Lyndon words. This result raises a natural question: are they the only ones? If no, what do the infinite smooth words that are also Lyndon words look like? In this paper, we give the answer, proving that the only infinite smooth Lyndon words are $m_{\{a<b\}}$, with $a,b$ even,  $m_{\{1<b\}}$ and $\Delta^{-1}_1(m_{\{1<b\}})$, with $b$ odd, where $m_\A$ is the minimal infinite smooth word with respect to the lexicographic order over a numerical alphabet $\A$ and $\Delta$ is the run-length encoding function.   
\end{abstract}


\section{Introduction}

Finite Lyndon words are easy to define: they are finite words smaller than any of their suffixes. They were first introduced by Chen, Fox and Lyndon in \cite{cfl1958} for the construction of bases of the lower central series for free groups. The authors proved that any finite word can be expressed as a unique non-increasing product (concatenation) of Lyndon words, which they called the Lyndon factorization.  Lyndon words were then studied by Duval \cite{jpd1983,jpd1988}. Among his results, is an algorithm that generates Lyndon words of bounded length for a finite alphabet and another one that computes the Lyndon factorization in linear time with respect to the length of the word.  Siromoney et al. \cite{smds1994}  later defined infinite Lyndon words as infinite words smaller than any of their suffixes, in order to introduce Lyndon factorization of infinite words. Lyndon words also appeared for instance in \cite{ml1983,mr1989,cr1993}. The Lyndon factorization gives nice properties about the structure of words. Since a few years, a wide literature is devoted to Lyndon words: \cite{bdl1997,hr2003,cr2005,gr2007,ps2003}. In particular, Melan\c con \cite{gm2000} studied Lyndon factorization of Sturmian infinite words and gave a formula that completely describes the Lyndon factorization of any characteristic Sturmian word.

The smooth infinite words over $\A=\{1,2\}$ form an infinite class {$\mathcal K$} of infinite words containing the well-known Kolakoski word $K$ \cite{wk1965} defined as one of the two fixed points of  the run-length encoding function $\Delta$,  that is
$$\Delta(K) = K = 2211212212211211221211212211211212212211212212 \cdots .$$  
They are characterized by the property that the orbit obtained by iterating  $\Delta$ is contained in $\{1,2\}^\omega$.  In  the early work of Dekking \cite{fmd1980}, there are some challenging conjectures on the structure of $K$ that still remain unsolved despite the efforts devoted to the study of patterns in $K$. For instance, we know from Carpi \cite{ac1994} that $K$ and more generally, any word in the infinite class $\mathcal K$ of smooth words over $\A=\{1,2\}$,  contain only a finite number of squares, implying by direct inspection that $K$ and any $w \in \mathcal K$ are cube-free. Weakley \cite{wk1965} showed that the  number of factors of length $n$  of ${\mathcal K}$  is  polynomially bounded.
In \cite{bl2003}, a connection was established between the palindromic complexity and the recurrence of $K$. Then Berth\'e, Brlek and Choquette \cite{bbc2005} studied smooth words over arbitrary alphabets and obtained a new characterization of the infinite Fibonacci word $F$.  Relevant work may also be found in \cite{bblp2003} and in \cite{bbc2005,jp2005}, where generalized Kolakoski words are studied for arbitrary alphabets. The authors investigated in \cite{bpm2007} the extremal infinite smooth words, that is the minimal and the maximal ones w.r.t. the lexicographic order, over the alphabets $\{1,2\}$ and $\{1,3\}$: a surprising link is established between $F$ and the minimal infinite smooth word over $\{1,3\}$.

More recently, Brlek, Jamet and Paquin  \cite{bjp2006} studied the extremal smooth words for any $2$-letter numerical alphabet and  they showed the existence of infinite smooth words that are also Lyndon words: the minimal smooth word over an even alphabet and the one over the alphabet $\{1,b\}$, with $b$ odd, are  Lyndon words. Then a natural question arises: are there other infinite smooth words that are also infinite Lyndon words? 

In this paper, we show that the only infinite smooth words that are also Lyndon words are related to the minimal smooth Lyndon words given in  \cite{bjp2006}. In order to prove it, we study the words over a 2-letter alphabet depending on the parity of the letters. The paper is organized as follows. In Section 2, we recall the basic definitions in combinatorics on words, we state the notation we will use next and we give useful known results.  Section 3 is devoted to the characterization of infinite smooth Lyndon words. It is divided in $4$ subsections. In Section 3.1, we study the case of an alphabet $\A=\{a<b\}$, with $a$ even and $b$ odd. We show that there is no infinite Lyndon words that is also smooth. In Section 3.2, we are interested in even alphabets. We show that only the minimal smooth word is a Lyndon word. Section 3.3 is devoted to odd alphabet. We prove that only $m_{\{1,b\}}$ $\Delta^{-1}_1(m_{\{1,b\}})$ are Lyndon word. Finally, Section 3.4 studies the words over an alphabet $\{a<b\}$, with $a$ odd and $b$ even. In this last case, we show that there is no infinite Lyndon word that is also smooth. 

This paper is an extended version of a paper presented in Prague (Czech Republic) during the 13th Prague Stringology Conference \cite{gp20082}. In this new version, we give the complete proofs. We also correct our main result, since Proposition 24 in \cite{gp20082} appeared to be false. That leads to the third form of smooth Lyndon words, namely $\Delta^{-1}_1(m_{\{1,b\}})$, with $b$ odd. 

\section{Preliminaries}

Throughout this paper,  $\A$ is a finite {\em alphabet} of {\em letters} equipped with a total order $<$.
A  {\em finite word} $w$ is a finite sequence of letters $w=w[0]w[1]\cdots w[n-1]$, where $w[i]\in \A$ denotes its $(i+1)$-th letter. Its length is $n$ and we write $|w|=n$. The set
of $n$-length words over $\A$ is denoted by $\A^n$. By
convention the {\it empty} word is denoted by $\varepsilon$ and its
length is $0.$ The free monoid generated by $\A$ is defined by
 $\A^* = \bigcup_{n\geq 0} \A^n$ and $\A^* \setminus \varepsilon$ is denoted $\A^+$.
 The set of {\it right infinite words}, also called {\it infinite words} for short, is denoted by
$\A^{\omega}$ and $\A^{\infty} = \A^*\cup
\A^{\omega}$. Adopting a consistent notation for finite words over the infinite alphabet $\mathbb{N}$, $\mathbb{N}^* = \bigcup_{n\geq 0} \mathbb{N}^n$ is the set of finite
words and  $\mathbb{N}^\omega$ is that of infinite ones. Given a word
$w \in \A^*,$  a {\sl factor} $f$ of $w$ is a word
$f \in \A^*$ satisfying
$$\exists  x,y \in \A^*, w=xfy .$$ If $x=\varepsilon$ (resp. $y=\varepsilon$ ) then $f$ is called a {\em prefix}
(resp. {\em suffix}). Note that by convention, the empty word is suffix and prefix of any word. A  {\em block} of length $k$ is a maximal factor of the particular form
$f =\alpha^k$, with $\alpha \in \A$.
The set of all factors of $w$, also called the {\it language} of $w$, is denoted by
$F(w)$ and those of length  $n$ is $F_n(w) = F(w) \cap
\A^n$. We denote by $\Pref(w)$ (resp. $\Suff(w)$)  the set of all prefixes (resp. suffixes) of
$w$.

Over an arbitrary 2-letter alphabet $\A= \{a, b\}$, there is
a usual length preserving morphism, the \textit{complementation},
defined by $ \overline{a}=b \; , \; \overline{b}=a,$ which
extends  to words as follows. The complement of $u = u[0]u[1]\cdots u[n-1] \in\A^n$ is the word $\overline{u}=
\overline{u[0]}\; \overline{u[1]} \cdots \overline{u[n-1]}$. The  {\it reversal} of $u$ is the word $\widetilde{u}= u[n-1]\cdots u[1]u[0].$

For $u, v \in \A ^*$, we write $u < v$ if and
only if $u$ is a proper prefix of $v$ or if there exists an integer $k$ such
that $u[i] = v[i]$ for $ 0\leq i\leq k-1$ and $u[k] < v[k]$. The relation
$\leq$ defined by $u\leq v$ if and only if $u=v$ or $u < v$, is
called the {\it lexicographic order}. That definition holds for $\A^{\omega}$.
Note that in general the complementation does not preserve  the lexicographic order. Indeed, when $u$ is not a proper prefix of $v$ then \begin{equation}\label{lexico}
 u > v \iff  \overline{u} < \overline{v}.
\end{equation}

\subsection{Lyndon words and factorizations}

A word $u \in \A ^*$ is a {\em Lyndon word} if  $u < v$ for all proper non-empty suffixes $v$ of $u$.
For instance, the word $11212$ is a Lyndon word while $12112$ is not since $112 < 12112$. 
A word of length 1 is clearly a Lyndon word. The set of Lyndon words is denoted by $\mathcal L.$

From the works of Chen, Fox and Lyndon, we have the following theorem.

\begin{theorem} \textnormal{\cite{cfl1958}}\label{finitelyndon} Any non empty finite word $w$ is uniquely expressed as a non increasing product of Lyndon words
 \begin{equation}\label{lyndon}
 w=\ell_0\ell_1\cdots\ell_n = \bigodot_{i=0}^n \ell_i,
 \;\textnormal{with}\quad \ell_i  \in {\mathcal L} \;\textnormal{ and    } \ell_0\geq \ell_1\geq \cdots \geq \ell_n .
 \end{equation}
\end{theorem}

Siromoney et al. \cite{smds1994} extended Theorem \ref{finitelyndon} to infinite words. The set $\mathcal L _\infty$ of {\it infinite Lyndon words} consists of infinite words smaller than any of their suffixes.

\begin{theorem} \textnormal{\cite{smds1994}} Any infinite word $w$ is uniquely expressed as a non increasing product of Lyndon words, finite or infinite, in one of the two following forms:
\begin{quote}\begin{itemize}
\item[{\rm i)}] either there exists an infinite sequence $(\ell_k)_{k \geq 0}$ of elements in ${\mathcal L}$  such that 
\begin{center}
\smallskip
$w=\ell_0\ell_1\ell_2\cdots$ and for all $\,k , \,\ell_k\geq \ell_{k+1};$
\end{center}
\item[{\rm ii)}] there exist  a finite sequence $\ell_0,\ldots,\ell_{m}  \,(m\geq 0)$ of elements in ${\mathcal L}$ and $\ell_{m+1}\in \mathcal L_\infty$ such that
\begin{center}
\smallskip
$w=\ell_0\ell_1\cdots\ell_m\ell_{m+1}$ and  $\ell_0 \geq \cdots \geq \ell_{m} >\ell_{m+1}. $
\end{center}
\end{itemize}
\end{quote}
\end{theorem}

Let us recall from (\cite{ml1983} Chapter 5.1) a useful property concerning Lyndon words.

\begin{lemma} \label{lyndonconc} Let $u,v \in \mathcal L$. Then $uv \in \mathcal L$ if and only if $u < v$.
\end{lemma}

A direct corollary of this lemma is:

\begin{corollary} \label{propb} Let $u, v \in \mathcal L$, with $u < v$. Then $\displaystyle uv^n, \, u^nv  \in \mathcal L,$  for all $n \geq 0$. 
\end{corollary}

\subsection{Run-length encoding}

The widely known {\it run-length encoding}  is used in many
applications as a method  for compressing data. For instance,  the
first step in the algorithm used for compressing the  data
transmitted by Fax machines consists of a run-length encoding of
each line of pixels.  Let $\A=\{a<b\}$ be an ordered alphabet.  Then every word  $w\in \A^*$ can
be uniquely written as a product of factors as follows:
$$w={a}^{i_0}{b}^{i_1}{a}^{i_2}\cdots \qquad \textnormal{or } \qquad w={b}^{i_0}{a}^{i_1}{b}^{i_2}\cdots  $$
with $i_k \geq 1$ for $k \geq 0$.
The operator giving the size of the blocks appearing in the coding is a function
$\Delta : \A^* \longrightarrow \mathbb{N}^*,$ defined by 
$\Delta(w)=i_0,i_1, i_2, \cdots $
which is easily extended to infinite words as
$\Delta :\A^{\omega} \longrightarrow \mathbb{N}^{\omega}.$

For instance, let $\A=\{1,3\}$ and $w = 13333133111$. Then
$$      w  = 1^1 3^4 1^1 3^2 1^3\quad \hbox{\rm and}\quad
        \Delta(w)  =  [1,4,1,2,3].$$
When $\Delta(w) \subseteq \{1,2,\cdots,9\}^*$, the punctuation and the parentheses are often omitted in order to
manipulate the more compact notation
$\Delta(w)=14123.$
This  example is a special case where the coding integers do not coincide
with the alphabet on which is encoded $w$, so that $\Delta$ can be viewed
as a partial function
$\Delta: \{1,3\}^* \longrightarrow \{1,2,3,4\}^*.$\\

From now on, we only consider 2-letter alphabets $\A=\{a<b\}$, with $a, b \in \N \setminus \{0\}$.\\

Recall from \cite{bl2003} that  $\Delta$ is not bijective since $\Delta(w)= \Delta(\overline{w}), $ but commutes with the reversal
$(\widetilde{\;\;\;\;} )$,  is stable under
complementation $(\,\raisebox{6pt}{$\overline{\;\;}$}\,)$ and preserves palindromicity.
Since   $\Delta$ is not bijective, pseudo-inverse functions 
$$\Delta_a^{-1},  \Delta_b^{-1} : \A^* \longrightarrow \A^*$$ 
are defined for 2-letter alphabets by
$$\Delta_\alpha^{-1}(u)=
\alpha^{u[0]}\overline{\alpha}^{u[1]}\alpha^{u[2]}\overline{\alpha}^{u[3]}
\cdots,\quad {\textnormal{for}}\,\,\, \alpha\in\{a,b\}.$$
Note that  the pseudo-inverse function $\Delta^{-1}$ also commutes with the mirror image, that is 
\begin{equation}
\widetilde{\Delta ^{-1}_\alpha(w)}=\Delta ^{-1}_\beta(\widetilde{w}), \label{miroir}
\end{equation} 
where $\beta = \alpha$ if $|w|$ is odd and $\beta=\overline{\alpha}$ if $|w|$ is even.

The operator $\Delta$ may be iterated, provided the process is stopped
when the coding alphabet changes or when the resulting word has length $1$.
\begin{example} \label{ex1}
{\rm Let $w=1333111333133311133313133311133313331113331$. The successive application of $\Delta$ gives :}\smallskip\\
$ \makebox[15pt]{ } \Delta^0(w)={\bf 1}333111333133311133313133311133313331113331;\\
 \makebox[15pt]{ }  \Delta^1(w)={\bf 1} 333133311133313331;\\
 \makebox[15pt]{ }  \Delta^2(w)={\bf 1} 31333131;\\
 \makebox[15pt]{ }  \Delta^3(w)={\bf 1} 113111;\\
 \makebox[15pt]{ }  \Delta^4(w)={\bf 3} 13;\\
 \makebox[15pt]{ }  \Delta^5(w)={\bf 1} 11;\\
 \makebox[15pt]{ }  \Delta^6(w)={\bf 3} .$
\end{example}

\subsection{Smooth words}

The set of {\it finite smooth words} over the alphabet $\A=\{a<b\}$ is defined by
$$\Delta^+_\A=\{w \in \A^*\,  | \, \exists n \in \N, \Delta^n(w)=\alpha \in \A \textnormal{ and }  \forall k \leq n, \Delta^k(w) \in \A^*\}.$$ 

The operator $\Delta$ extends to infinite words (see \cite{bl2003}). Define the set of {\it infinite smooth words} over $\A=\{a<b\}$ by
\[{\mathcal K}_{\A}= \{ w \in \A^\omega\;|\;
\forall k \in \N, \Delta^k(w) \in \A^{\omega}\}.\]

The well-known \cite{wk1965} {\it Kolakoski word} denoted $K$ is defined as the fix-point starting with the letter $2$ of the operator $\Delta$ over the alphabet $\{1,2\}$:
$$K=2211212212211211221211212211211212212211\cdots$$ 

More generally, the operator $\Delta$ has two fix-points in $\mathcal{K}_\A$, namely
\[\Delta(K_{(a,b)})=K_{(a,b)} \quad \textnormal{and }\quad \Delta(K_{(b,a)}) = K_{(b,a)} ,\]
where $K_{(a,b)}$ is the generalized Kolakoski word \cite{jp2005} over the alphabet $\{a,b\}$ starting with the letter $a$. 

\begin{example} The Kolakoski word  over  $\A=\{1,2\}$ starting with the letter $2$ is $K~=~K_{(2,1)}$. We also have $K_{(2,3)}= 22332223332233223332\cdots$   and    $K_{(3,1)}=33311133313133311133\cdots$.
\end{example}

A bijection $\Phi : \mathcal{K}_{\A} \longrightarrow \A ^\omega$ is defined by
$$\Phi(w)=\Delta^0(w)[0]\Delta^1(w)[0]\Delta^2(w)[0]\cdots =\prod_{i\geq 0 }\Delta^i (w)[0]$$
and its inverse is defined as follows. Let $u \in \A^k$, then $\Phi ^{-1}(u)=w_k,$
where 
\begin{displaymath}
w_n = \left \{ 
\begin{array}{ll}
u[k-1], & \textnormal{if } n = 1;\\
\Delta^{-1}_{u[k -n]}(w_{n-1}), & \textnormal{if } 2 \leq n \leq k.
\end{array} \right.
\end{displaymath}
Then for $k=\infty$, $\Phi^{-1}(u)=\lim_{k \rightarrow \infty} w_k=\lim_{k\rightarrow \infty} \Phi^{-1}(u[0..k-1])$.

\begin{remark} With respect to the usual topology defined by 
$$d((u_n)_{n\geq 0}, (v_n)_{n \geq 0}):= 2^{-\min\{j \in \N, u_j\neq v_j \}},$$ 
the previous limit exists since $w_k$ is prefix of $w_{k+1}$ for all $k \in \N$.
\end{remark}

\begin{example} For the word $w=1333111333133311133313133311133313331113331$ of the Example \ref{ex1}, $
\Phi(w)= {\bf 1111313}$. 
\end{example}

Note that since $\Phi$ is a bijection between the infinite smooth words and the infinite words,  the set of  infinite smooth words is infinite. Moreover, given a prefix of $\Phi(w)$ with $w$ a smooth word, we can construct a prefix of $w$ as in the following example.

\begin{example} \label{exForce} Let $p=1221$ be a prefix of $\Phi(w)$, with $w \in \{1,2\}^\omega$ an infinite smooth word. Then we compute from bottom to top, using the $\Delta^{-1}$ operator:\\
\indent $\Delta^0(w)=11221221 \cdots \\
\indent \Delta^1(w)=2212 \cdots \\
\indent \Delta^2(w)=21 \cdots \\
\indent \Delta^3(w)=1\cdots$ \\
Note that in $\Delta^2(w)$, the letter $1$ is obtained by deduction, since $\Delta^3(w)$ indicates that the first block of letters of $\Delta^2(w)$ has length $1$. The last written letter of every line is deduced by a similar argument.
\end{example}

We recall from \cite{bpm2007} the useful {\em right derivative}   $D_r: \A ^*  \rightarrow \N^*$ defined by \begin{displaymath}
D_r(w)=  \left \{ \begin{array}{ll}
\varepsilon & \textrm{if $\Delta(w)=\alpha$, $\alpha <b$ or $w=\varepsilon$,}\\
\Delta(w) & \textrm{if $\Delta(w)=xb$,}\\
x & \textrm{if $\Delta(w)=x\alpha$, $\alpha <b$,}
\end{array}  \right . 
\end{displaymath}
where $\alpha \in \N$ and $x \in \A^*$.
A word $w$ is {\it r-smooth} (also called a  {\it smooth prefix}) if $\forall k \geq 0, \; D_r^{k}(w) \in \A ^*$. In other words, if a word $w$ is {\it r-smooth}, then it is a prefix of at least one infinite smooth word (see \cite{gp2008} for more details). 

\begin{example} Let $w=112112212$. Then $\Delta(w)=212211$, $\Delta^2(w)=1122$, $\Delta^3(w)=22$, $\Delta^4(w)=2$ and $D_r(w)=21221$, $D^2_r(w)=112$, $D_r^3(w)=2$. 
\end{example}

Similarly, the operator $D$ is defined over the alphabet $\{a<b\}$ by
\begin{displaymath}
D(w)=  \left \{ \begin{array}{ll}
\varepsilon & \textrm{if $\Delta(w)<b$ or $w=\varepsilon$,}\\
\Delta(w) & \textrm{if $\Delta(w)=bxb$ or $\Delta(w)=b$,}\\
bx & \textrm{if $\Delta(w)=bxu$,} \\
xb & \textrm{if $\Delta(w)=uxb$,} \\
x & \textrm{if $\Delta(w)=uxv$,}
\end{array}  \right .
\end{displaymath}
where  $u$ and $v$ are blocks of length $< b$. 
A finite word is called  a {\it smooth factor} (also called a $C^\infty$-word in \cite{ac1994,ybh1997,ybh2008,wdk1989})  if there exists $k\in \N$ such that $D^k(w)=\varepsilon$ and  $\forall j <k, \; D^{j}(w)\in \A ^*$.

\subsection{Known results}

The {\it minimal} (resp. the {\it maximal}) {\it infinite smooth word} over the alphabet $\A$ is the smallest (resp. biggest) infinite smooth word, with respect to the lexicographic order. It is denoted by $m_\A$ (resp. $M_\A$). 

An alphabet $\A=\{a<b\}$ is called an {\it odd alphabet} (resp. {\it even alphabet}) if both $a$ and $b$ are odd (resp. even). The extremal smooth words satisfy the following properties established in a previous paper. 

\begin{proposition}  \textnormal{\cite{bjp2006}} \label{motMin} Let $\A=\{a<b\}$. Then the following properties hold.
\begin{itemize}
\item [\rm i)] If $a$ and $b$ are both even we have :\\
$\Phi(M_{\{a,b\}})= b^\omega$;  \quad $\Phi(m_{\{a,b\}})= ab^\omega$; \quad and \quad
$m_{\{a,b\}} \in \mathcal L_\infty$.
\item [\rm ii)] If $a$ and $b$ are both odd we have :\\
$\Phi(M_{\{a,b\}})= (ba)^\omega$;  \quad $\Phi(m_{\{a,b\}})= (ab)^\omega$; \quad and \quad
$m_{\{a,b\}} \in \mathcal L_\infty \iff a=1$.
\end{itemize}
\end{proposition}

Let us recall useful results about smooth words. 

\begin{lemma} \textnormal{\cite{bbc2005}}\label{glemma} Let $u$, $v$ be finite smooth words. If there exists an index $m$ 
such that, for all $i$, $0 \leq i \leq m$, the last letter of $\Delta ^i (u)$ differs from the first letter 
of $\Delta ^i (v)$, and $\Delta ^i (u)\neq 1$, $\Delta ^i(v)\neq 1$, then
\smallskip\\ 
 \makebox[25pt]{\hfill\rm i)} $\Phi(uv)=\Phi(u)[0..m] \cdot \Phi \circ \Delta ^{m+1}(uv)$;\\
 \makebox[25pt]{\hfill\rm ii)} $\Delta ^i (uv) = \Delta ^i (u)\Delta ^i  (v)$. 
\end{lemma}

The following properties follow immediately from the definitions. For more details, the reader is referred to \cite{gp2008}.
 
Recall from \cite{bdlv2006} that in the case of the alphabet $\A=\{1,2\}$, every finite smooth word $w \in \Delta^+_\A$ can be easily extended to the right in a smooth word by means of the function $\Phi$ as 
$\forall u \in \A^\infty, \, w \in \Pref(\Phi^{-1}(\Phi(w)\cdot u)).$ Its generalization to arbitrary alphabets is immediate (see \cite{gp2008}).

\begin{proposition} \textnormal{\cite{gp2008}} \label{propPhi} Let $\A=\{a<b\}$. Then the following properties hold.
\begin{itemize}
\item [\rm i)] Any smooth prefix can be arbitrarily right extended to an infinite smooth word. 
\item [\rm ii)] Let $u=\Phi(w)$, with $w\in \A^\omega$ an infinite smooth word. If $u=u'u''$, then $\Phi^{-1}(u')$ is prefix of $w$. 
\end{itemize}
\end{proposition}

\section{Characterization of infinite smooth Lyndon words}

In this section, we prove our main result: the only infinite smooth words that are also infinite Lyndon words are $m_{\{a<b\}}$, with $a, b$ even, $m_{\{1<b\}}$ and $\Delta^{-1}_1(m_{\{1<b\}})$, with $b$ odd.

In order to prove it, we study the four possible combinations of the parity of the letters $a$ and $b$. For each case, we fix a length $n$ and we then consider all the possible words $p$ of length $\leq n$ such that $\Phi^{-1}(p)$ is prefix of an infinite smooth word $w$. We suppose that $w$ is also a Lyndon word. Then for each word $p$, either  we prove that $\Phi^{-1}(p)$ can not be a prefix of a Lyndon word by showing the existence of a smaller suffix, or we describe an infinite smooth Lyndon word having $\Phi^{-1}(p)$ as prefix. For each case, the different values of $p$ are illustrated in a tree. 

\begin{lemma} \label{15} Let $p\in \{a<b\}^*$ be such that $\Phi^{-1}(p)$ is prefix of an infinite smooth Lyndon word $w$. Then $p[0]=a$. 
\end{lemma}

\Proof Follows from the equality $p[0]=w[0]$ and since a Lyndon word $w$ must start by the smallest letter, namely $a$. \QED

\noindent Lemma \ref{15} will be used in this section to exclude the cases numbered~(0) in the proofs. 

\subsection{Over $\A$ with $a$ even and $b$ odd}

In this section, we prove the following result.

\begin{theorem} Over the alphabet $\{a<b\}$, with $a$ even and $b$ odd, there is no infinite smooth word that is also a Lyndon word. 
\end{theorem}

\begin{figure}[h]
\begin{center}
\includegraphics[width=4.5cm]{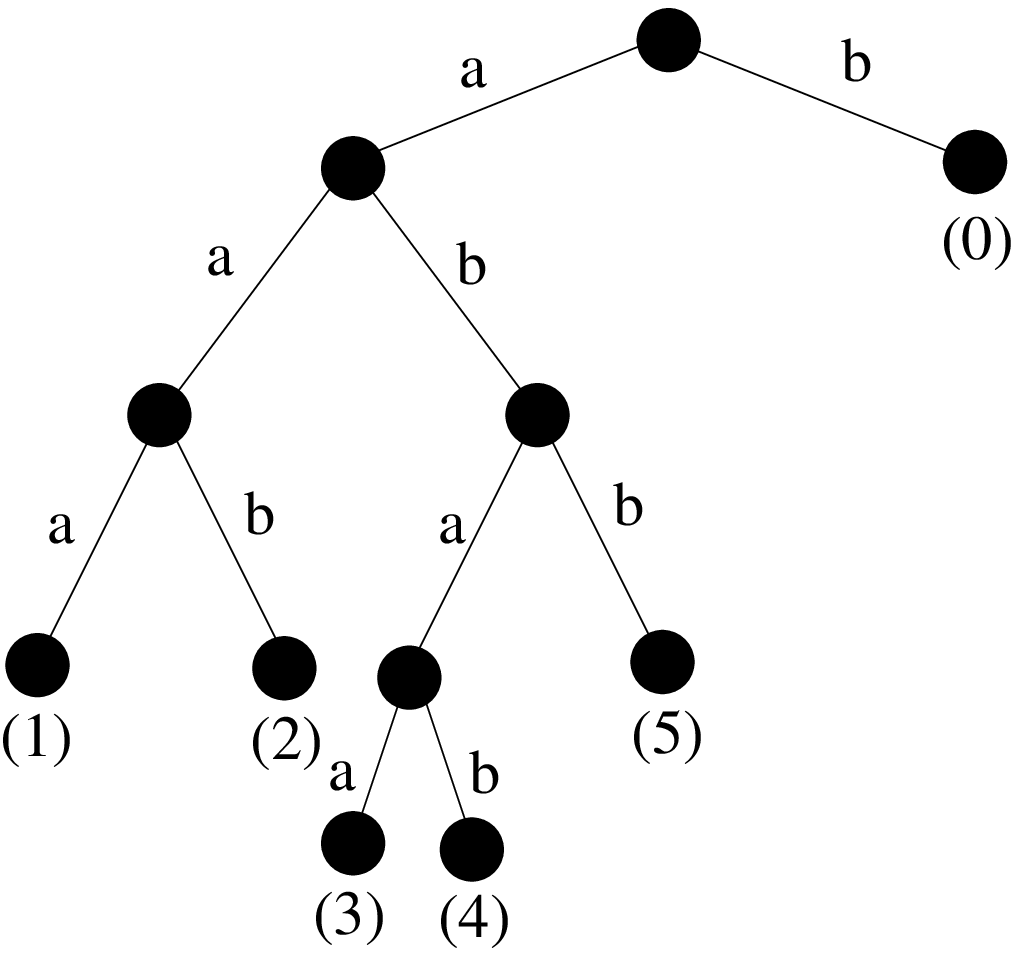}
\caption{Possible cases for an even-odd alphabet}\label{figpi}
\end{center}
\end{figure}

\Proof  Figure \ref{figpi} illustrates the $6$ possible cases to consider, using a tree. The leaves correspond to the first letter of $\Phi(w)$ that leads to a contradiction: the prefix $\Phi^{-1}(p)$ obtained can not be the prefix of an infinite Lyndon word. We will prove it by showing that there exists a factor $f$ of $w$ not prefix of $\Phi^{-1}(p)$ such that $f < w$. For clarity issues, the first letter of $f$ is underlined in $w$.  

\begin{enumerate}
\item [Case (1)]  If $p=aaa$, then \\
\indent $\displaystyle \Delta^0(w)= (a^ab^a)^\frac{a}{2}(\underline a^b b^b)^\frac{a}{2}\cdots\\
\indent \Delta^1(w)= a^a b^a \cdots\\
\indent  \Delta^2(w)=aa \cdots$\\
Since $w$ has the prefix $a^ab^a$ and the factor $f=a^b$, it can not be a Lyndon word. 

\begin{remark}\label{remxx} Since the smallest letter $a$ of the alphabet is even, $2 \leq a \leq p[3]$. That allows us to assume  that $\Delta^2(w)$ starts with a block of length at least $2$. This argument holds for $\Delta^i(w)$, $i \geq 0$, and will be used for almost all the cases considered in this paper. 
\end{remark}

\begin{remark} In the previous case, we construct $\Delta^0(w)$ from $\Delta^2(w)$, applying $\Delta^{-1}$ twice. We will always proceed this way.
\end{remark}

\item [Case (2)] If $p=aab$, then\\
\indent $\displaystyle \Delta^0(w)=(a^ab^a)^\frac{b-1}{2}a^a(b^b\underline a^b)^\frac{b-1}{2}b^b\cdots \\
\indent \Delta^1(w)=a^bb^b\cdots \\
\indent \Delta^2(w)=bb\cdots$\\\
$w$ has the factor $f=a^b$ smaller than its prefix $a^ab^a$. 

\item [Case (3)] If $p=abaa$, then\\
\indent $\displaystyle \Delta^0(w)=((a^bb^b)^\frac{a}{2}(a^ab^a)^\frac{a}{2}) ^\frac{a}{2}((\underline a^bb^b)^\frac{b-1}{2}a^b(b^aa^a)^\frac{b-1}{2}b^a)^\frac{a}{2}\cdots \\
\indent \Delta^1(w)=(b^aa^a)^\frac{a}{2}(b^ba^b)^\frac{a}{2}\cdots\\
\indent \Delta^2(w)=a^ab^a\cdots\\
\indent \Delta^3(w)=aa\cdots$\\
$w$ has the factor $f=(a^bb^b)^\frac{b-1}{2}a^b$.

\item [Case (4)] If $p=abab$, then \\
\indent $\displaystyle \Delta^0(w)=((a^bb^b)^\frac{a}{2}(a^ab^a)^\frac{a}{2})^\frac{b-1}{2}(a^b b^b)^\frac{a}{2}((a^ab^a)^\frac{b-1}{2}a^a(b^b a^b)^\frac{b-1}{2}b^b)^\frac{b-1}{2}(a^a b^a)^\frac{b-1}{2}a^a (b^b\underline a^b)^{\frac{a}{2}} b^aa\cdots \\
\indent \Delta^1(w)=(b^aa^a)^\frac{b-1}{2}b^a(a^bb^b)^\frac{b-1}{2}a^bb^aa \cdots \\
\indent \Delta^2(w)=a^bb^ba\cdots \\
\indent \Delta^3(w)=bb \cdots $\\
$w$ has the smaller factor $f=a^bb^a a$.

\begin{remark} \label{remx} In the computation of $\Phi^{-1}(abab)$, we deduce the last letter of each line, as we did in Example \ref{exForce}. This deduction will be used in further cases.
\end{remark}

\item [Case (5)] If $p=abb$, then \\
\indent $\displaystyle \Delta^0(w)=(a^bb^b)^\frac{b-1}{2}\underline a^b(b^aa^a)^\frac{b-1}{2} b^a\cdots \\
\indent \Delta^1(w)=b^ba^b \cdots \\
\indent \Delta^2(w)=   bb \cdots $\\
$w$ has the factor $f=a^bb^aa^a $. 
\end{enumerate}

We conclude using Proposition \ref{propPhi} ii). \QED

\subsection{Over an even alphabet}

Let us now consider the case of an alphabet $\A$ having even letters. 

\begin{theorem} Over the alphabet $\{a<b\}$, with $a$ and $b$ even, the only smooth word that is also an infinite Lyndon word is $m_{\{a<b\}}$. 
\end{theorem}

\begin{figure}[h]
\begin{center}
\includegraphics[width=4cm]{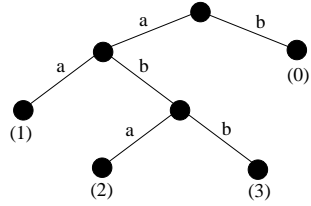}
\caption{Possible cases for an even alphabet}\label{figp}
\end{center}
\end{figure}

\Proof We proceed similarly as in the previous section.  The $4$ possibilities are illustrated in Figure \ref{figp}. We will first prove that Cases (1) and (2) are impossible and then, we will show in Case (3) that $abb$ is prefix of an infinite smooth Lyndon word $w$ if and only if $w=\Phi^{-1}(ab^\omega)=m_{\{a,b\}}$.

\begin{enumerate}
\item [Case (1)] If $p=aax$, with $x\in \A$, then\\
\indent $\displaystyle \Delta^0(w)=(a^ab^a)^\frac{x}{2}(\underline a^bb^b)^\frac{x}{2}\cdots \\
\indent \Delta^1(w)=a^xb^x\cdots \\
\indent \Delta^2(w)=xx\cdots $\\
$w$ has the factor $f=a^b$.

\item [Case (2)] If $p=abax$, with $x\in \A$, then\\
\indent $\displaystyle \Delta^0(w)=((a^bb^b)^\frac{a}{2}(a^ab^a)^\frac{a}{2})^\frac{x}{2}((\underline a^bb^b )^\frac{b}{2}(a^ab^a)^\frac{b}{2})^\frac{x}{2} \cdots \\
\indent \Delta^1(w)=(b^aa^a)^\frac{x}{2}(b^ba^b)^\frac{x}{2}\cdots \\
\indent \Delta^2(w)=a^xb^x \cdots\\
\indent \Delta^3(w)=xx\cdots $\\
$w$ has the factor $f=(a^bb^b)^\frac{b}{2}$.

\item [Case (3)] Recall that  the minimal infinite smooth word $m_{\{a,b\}}=\Phi^{-1}(ab^\omega)$ is a Lyndon word. Let us show that it is the only smooth word that is also a Lyndon word. In order to prove it,  let us suppose that  we can write $p=ab^kay$, with $k \geq 2$ maximal, since Case (2) already excluded the possibility $k=1$,  and $y \in \A^\omega$. Let us first compute $u=\Phi^{-1}(bbay)$, with $y[0]=x \in \A$. We get\\
 \indent  $\Delta^0(u)=((b^ba^b)^\frac{a}{2}(b^aa^a)^\frac{a}{2})^\frac{x}{2} ((b^ba^b)^\frac{b}{2}(b^aa^a)^\frac{b}{2})^\frac{x}{2}\cdots \\
 \indent \Delta^{1}(u)=(b^aa^a)^\frac{x}{2}(b^ba^b)^\frac{x}{2} \cdots \\
 \indent \Delta^{2}(u)=a^xb^x \cdots\\
 \indent \Delta^3(u)=xx \cdots $\\

Since $a$ and $b$ are even and using Lemma \ref{glemma},  $\Phi^{-1}(b^kay)$ can be then written as 
$$(v_1^{\frac{a}{2}}v_2^{\frac{a}{2}})^{\frac{x}{2}}(v_1^{\frac{b}{2}}v_2^{\frac{b}{2}})^{\frac{x}{2}}s,$$
with $v_1=\Delta^{-(k-2)}_b(b^ba^b)$, $v_2=\Delta_b^{-(k-2)}(b^aa^a)$, $s \in \A^\omega$ and $x \in \A$. 

Moreover, since $\Phi^{-1}(b^\omega)$ is the maximal smooth word and since $\Phi^{-1}(b^k)$ (resp. $\Phi^{-1}(b^{k-1}a)$) is prefix of $v_1$ (resp. $v_2$), we have that $v_1>v_2$ and $v_1$ is not prefix of $v_2$. Furthermore for $k\geq 2$, using Equation (\ref{lexico}) we get
$$\Delta^{-1}_b(v_1) > \Delta^{-1}_b(v_2) \iff \Delta^{-1}_a(v_1) < \Delta^{-1}_a(v_2),$$ 
that implies
$$ \Delta^{-1}_a(v_1^{\frac{b}{2}})< \Delta^{-1}_a(v_1^{\frac{a}{2}}v_2^{\frac{a}{2}}).$$
Thus $\Phi^{-1}(ab^kay)$ has the prefix $\Delta^{-1}_a(v_1^{\frac{a}{2}}v_2^{\frac{a}{2}})$ and the smaller factor not prefix $\Delta^{-1}_a(v_1^{\frac{b}{2}})$. Consequently, $\Phi^{-1}(ab^kay)$ is not a Lyndon word. Finally, $pp'=abbp'$ is such that $\Phi^{-1}(pp')$ is a smooth Lyndon word if and only if $p'=b^\omega$. 
\end{enumerate}
The only smooth Lyndon word over  an even 2-letter alphabet is the minimal smooth word $m_\A$ with $\Phi(m_\A)=ab^\omega$.  \QED

\subsection{Over an odd alphabet}

In this section, we prove the following result.

\begin{theorem} \label{thmii} Over the alphabet $\{a<b\}$, with $a$ and $b$ odd, there exist $2$ infinite smooth Lyndon words if and only if $a=1$. More precisely, the smooth Lyndon word are the minimal smooth word $m_{\{1<b\}}$ and $\Delta^{-1}_1(m_{\{1,b\}})$.
\end{theorem}

Before proving Theorem \ref{thmii}, some results are required. 

\begin{lemma} \label{lemdelta} Let $\A=\{a<b\}$ be an odd alphabet. Let $w,w'$ be two factors of a smooth word such that $w< w'$ and $w=xay$,  $w'= xby'$, with $x, y, y' \in \A^*$. Then if $|x|$ is even (resp. odd),
$$\Delta^{-1}_{\alpha}(w)<\Delta^{-1}_{\alpha}(w') \iff \overline \alpha < \alpha \textnormal{ (resp. } \alpha < \overline \alpha),$$
with $\alpha \in \A$ and $\overline \alpha$ its complement. 
\end{lemma}

\Proof Assume $|x|$ even.  Then $\Delta^{-1}_\alpha(x)$ ends by $\overline \alpha$. By direct computation, we have the following equalities:
\begin{eqnarray*}
&\Delta^{-1}_{ \alpha}(w)= \Delta^{-1}_{ \alpha}( xay)= \Delta^{-1}_{ \alpha}( x) \Delta^{-1}_{ \alpha}(a) \Delta^{-1}_{\overline \alpha}(y)= \Delta^{-1}_{ \alpha}( x){ \alpha}^a \Delta^{-1}_{\overline \alpha} (y)&\\
& \textnormal{and}&\\
&\Delta^{-1}_{ \alpha}(w')= \Delta^{-1}_{ \alpha}( xby')= \Delta^{-1}_{ \alpha}( x) \Delta^{-1}_{ \alpha}(b) \Delta^{-1}_{\overline \alpha}(y')=\Delta^{-1}_{ \alpha}( x){ \alpha}^b \Delta^{-1}_{\overline \alpha}(y').&
\end{eqnarray*}

Then  $\Delta^{-1}_\alpha(w)<\Delta^{-1}_\alpha (w')$ if and only if $\Delta^{-1}_{\overline \alpha}(y)[0] < \alpha$.  We conclude since $\Delta^{-1}_{\overline \alpha}(y)[0]=\overline \alpha$. A similar argument holds for $|x|$ odd. \QED 

Now, let us prove $2$ subcases of Theorem \ref{thmii}: $a\neq 1$ (Theorem \ref{thmii1}) and $a=1$ (Theorem~\ref{thmOuf}).

\begin{theorem} \label{thmii1} Over the alphabet $\{a<b\}$, with $a,b$ odd and $a\neq 1$, there is no infinite smooth word that is also a Lyndon word.
\end{theorem}

\begin{figure}[h]
\begin{center}
\includegraphics[width=2.5cm]{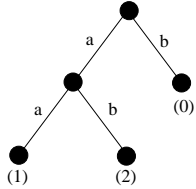}
\caption{Possible cases for an odd alphabet, with $a\neq 1$}\label{figii1}
\end{center}
\end{figure}

\Proof As in Sections $3.1$ and $3.2$, we proceed by inspection of the different possible prefixes of $\Phi(w)$ (see Figure \ref{figii1}) for an infinite smooth word $w$. 

\begin{enumerate}
\item [Case (1)] If $p=aax$, with $x \in \A$, then \\
\indent $\displaystyle \Delta^0(w)=(a^a b^a)^\frac{x-1}{2}a^a(b^b\underline a^b)^\frac{x-1}{2}b^b(a^a b^a)^\frac{x-1}{2}a^a\cdots \\
\indent \Delta^1(w)=a^xb^x a^x\cdots\\
\indent \Delta^2(w)=xxx \cdots$\\
Since $1\leq a \leq x$ and $x$ is odd, $x\geq 3$ and then, $\frac{x-1}{2}\geq 1$. Thus, $w$ has the factor $f=a^b$. 

\begin{remark} In the same way as in  Remark \ref{remxx}, we can suppose that $\Delta^2(w)$ starts by a block of length at least $3$. 
\end{remark}

\item [Case (2)] If $p=abx$, with $x \in \A$, then \\
\indent $\displaystyle \Delta^0(w)= (a^bb^b)^\frac{x-1}{2}\underline a^b(b^aa^a)^\frac{x-1}{2}b^a(a^bb^b)^\frac{x-1}{2}a^b\cdots\\
\indent \Delta^1(w)=b^x a^x b^x\cdots \\
\indent \Delta^2(w)=xxx\cdots$\\
$w$ has the factor $f=a^bb^aa^a$. 
\end{enumerate}

Since in the $3$ cases, it is possible to find a factor smaller than the prefix, we conclude that there is no smooth Lyndon word over an odd alphabet $\A$, with $a\neq1$. \QED

\begin{lemma} \label{lemPos} Let $w$ be a smooth word over an odd alphabet $\{a<b\}$, starting by the letter $a$ (resp. $b$). Then, a block of $a$'s always starts in an even (resp. odd) position and a block of $b$'s always starts in an odd (resp. even) position. 
\end{lemma}

\Proof It follows from the parity of the letters of the alphabet: every block has odd length, thus between $2$ blocks there is an even number of positions. \QED 

\begin{lemma} \label{blocb} In the smooth word $w=\Phi^{-1}((1b)^\omega)$, the blocks of $b's$ have all length $1$. \QED
\end{lemma}

\Proof It is already proved for $b=3$ in \cite{bpm2007} and follows from the parity of the letters and the alternance of the letters in $\Phi(w)$. \QED

We know from \cite{bjp2006} that the minimal smooth word over an odd 2-letter alphabet $\{1<b\}$ is a Lyndon word. The next proposition shows that this is the only infinite Lyndon word $w$ over the alphabet $\{1<b\}$ such that $\Phi(w)$ starts by $1b$. 

\begin{proposition} \label{propOuf1} Over the alphabet $\{1<b\}$, with $b$ odd, the only infinite smooth Lyndon word $w$ such that $\Phi(w)$ starts by $1b$ is $m_{\{1<b\}}$.
\end{proposition}

\begin{figure}[h]
\begin{center}
\psfrag{(1b)}{$(1b)^k$}
\includegraphics[width=4cm]{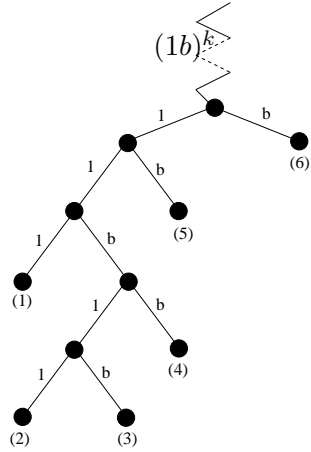}
\caption{Possible prefixes of $p$ starting by $(1b)^k$, $k \geq 1$}\label{fig1b}
\end{center}
\end{figure}

\Proof Recall first that $\Phi(m_{\{1,b\}})=(1b)^\omega$. Let us consider all possible prefixes $p$ such that $\Phi^{-1}(p)$ is prefix of an infinite smooth Lyndon word. Since we suppose that $p$ starts by $1b$, let us rewrite $p$ as $(1b)^ks$, where $k$ is maximal and $s \in \A^\omega$. We then proceed by inspection of the different possibilities (see Figure \ref{fig1b}).

\begin{enumerate}

\item [Case (1)] If $p=(1b)^k111$, let us consider the prefix  $\Phi^{-1}(1b111)$ of $u$.\\
 \indent  $\Delta^0(u)=1^b((b1)^{\frac{b-1}{2}}b(1^bb^b)^{\frac{b-1}{2}}\underline 1^b)^{\frac{b-1}{2}}(b1)^{\frac{b-1}{2}}b1^bb1\cdots \\
 \indent   \Delta^1(u)=b(1^bb^b)^{\frac{b-1}{2}}1^bb1\cdots\\
  \indent  \Delta^2(u)=1 b^b1\cdots \\
   \indent \Delta^3(u)= 1b \cdots\\
    \indent  \Delta^4(u)=1\cdots$\\
$u$ has the  factor $f=1^b(b1)^{\frac{b-1}{2}}b1^bb1$. Using Lemma \ref{lemdelta}  $2(k-1)$ times, Lemmas \ref{lemPos} and \ref{blocb}, we conclude that the prefix $p$ does not describe a prefix of a Lyndon  word.

\item [Case (2)] If $p=(1b)^k11b11$, let us consider the prefix $\Phi^{-1}(1b11b11)$ of $u$.\\
 \indent  $\displaystyle \Delta^0(u)=(1^b(b1)^{\frac{b-1}{2}}b)^{\frac{b-1}{2}}1^b(((b1)^{\frac{b-1}{2}}b(1^bb^b)^{\frac{b-1}{2}}1^b)^{\frac{b-1}{2}}(b1)^{\frac{b-1}{2}}b1^b)^{\frac{b-1}{2}}(b1^b)^{\frac{b-1}{2}}b((1^bb^b)^{\frac{b-1}{2}}\underline 1^bb)^{\frac{b-1}{2}}\\ 
 \indent \qquad \qquad (1^bb^b)^{\frac{b-1}{2}}1^b(b1^b)^{\frac{b-1}{2}}b(1^bb^b)^{\frac{b-1}{2}}1^bb1 \cdots \\
 \indent \Delta^1(s)=  (b1^b)^{\frac{b-1}{2}}b((1^bb^b)^{\frac{b-1}{2}}1^bb)^{\frac{b-1}{2}}(1b)^{\frac{b-1}{2}}1(b^b1)^{\frac{b-1}{2}}b^b(1b)^{\frac{b-1}{2}}1b^b1\cdots\\
  \indent  \Delta^2(u)=(1b)^{\frac{b-1}{2}}1(b^b1)^{\frac{b-1}{2}}1^b(b1)^{\frac{b-1}{2}}b1^bb\cdots \\
   \indent \Delta^3(u)=1^b(b1)^{\frac{b-1}{2}}b1^bb\cdots \\
    \indent  \Delta^4(u)=b1^bb\cdots \\
    \indent \Delta^5(u)=1b \cdots\\
    \indent \Delta^6(u)=1 \cdots$\\
$u$ has the factor $f=1^bb(1^bb^b)^{\frac{b-1}{2}}$. Applying Lemma \ref{lemdelta} $2(k-1)$ times and using Lemmas \ref{lemPos} and \ref{blocb}, we conclude.

\item [Case (3)] If $p=(1b)^k11b1b$, let us consider the prefix $\Phi^{-1}(11b1b)$ of $u$.\\
 \indent $\displaystyle \Delta^0(u)=((1b)^{\frac{b-1}{2}}1b^b)^{\frac{b-1}{2}}(1b)^{\frac{b-1}{2}}1(b^b1)^{\frac{b-1}{2}}1^b(b1)^{\frac{b-1}{2}}b1^bb \cdots \\
 \indent  \Delta^1(u)= (1^bb)^{\frac{b-1}{2}}1^b(b1)^{\frac{b-1}{2}}b1^bb\cdots \\
  \indent \Delta^2(u)=  (b1)^{\frac{b-1}{2}}b1^bb\cdots \\
   \indent \Delta^3(u)=1^bb\cdots \\
    \indent \Delta^4(u)= b \cdots $\\
$u$ has the factor $f=1^b$. Applying Lemma \ref{lemdelta} $(2k)$ times and using Lemmas \ref{lemPos} and \ref{blocb}, we conclude.

\item [Case (4)] If $p=(1b)^k11bby$,  let  us consider the prefix $\Phi^{-1}(1b11bb)$ of $u$.\\
 \indent  $\displaystyle \Delta^0(u)=((1^b(b1)^{\frac{b-1}{2}}b)^{\frac{b-1}{2}}1^b(((b1)^{\frac{b-1}{2}}b(1^bb^b)^{\frac{b-1}{2}}\underline 1^b)^{\frac{b-1}{2}}(b1)^{\frac{b-1}{2}}b(1^bb)^{\frac{b-1}{2}}1^b)^{\frac{b-1}{2}}\cdots\\
 \indent  \Delta^1(u)= ((b1^b)^{\frac{b-1}{2}}b((1^bb^b)^{\frac{b-1}{2}}1^b(b1)^{\frac{b-1}{2}}b)^{\frac{b-1}{2}}(1^bb^b)^{\frac{b-1}{2}}1^b)^{\frac{b-1}{2}}(b1^b)^{\frac{b-1}{2}}b(1^bb^b)^{\frac{b-1}{2}}1^bb1^bb\cdots \\
  \indent \Delta^2(u)= ((1b)^{\frac{b-1}{2}}1(b^b1^b)^{\frac{b-1}{2}}b^b)^{\frac{b-1}{2}}(1b)^{\frac{b-1}{2}}1b^b1b\cdots \\
   \indent  \Delta^3(u)=(1^bb^b)^{\frac{b-1}{2}}1^bb1 \cdots \\
    \indent  \Delta^4(u)=b^b1\cdots \\
    \indent \Delta^5(u)=b\cdots $\\
$u$ has the factor $f=1^b(b1)^{\frac{b-1}{2}}b(1^bb)^{\frac{b-1}{2}}$. Applying Lemma \ref{lemdelta} $2(k-1)$ times and using Lemmas \ref{lemPos} and \ref{blocb}, we conclude. 

\item [Case (5)] If $p=(1b)^k(1b)$, it is a prefix of $m_{\{1,b\}}$, which is a smooth Lyndon word.

\begin{figure}[b]
\begin{center}
\psfrag{1}{$1$}
\psfrag{v}{$v$}
\psfrag{2k}{$2k$}
\psfrag{b=(b-1)+1}{$b=(b-1)+1$}
\psfrag{b}{$b$}
\psfrag{vdot}{$\vdots$}
\psfrag{b(b-1)}{$b^{b-1}$}
\includegraphics[width=5cm]{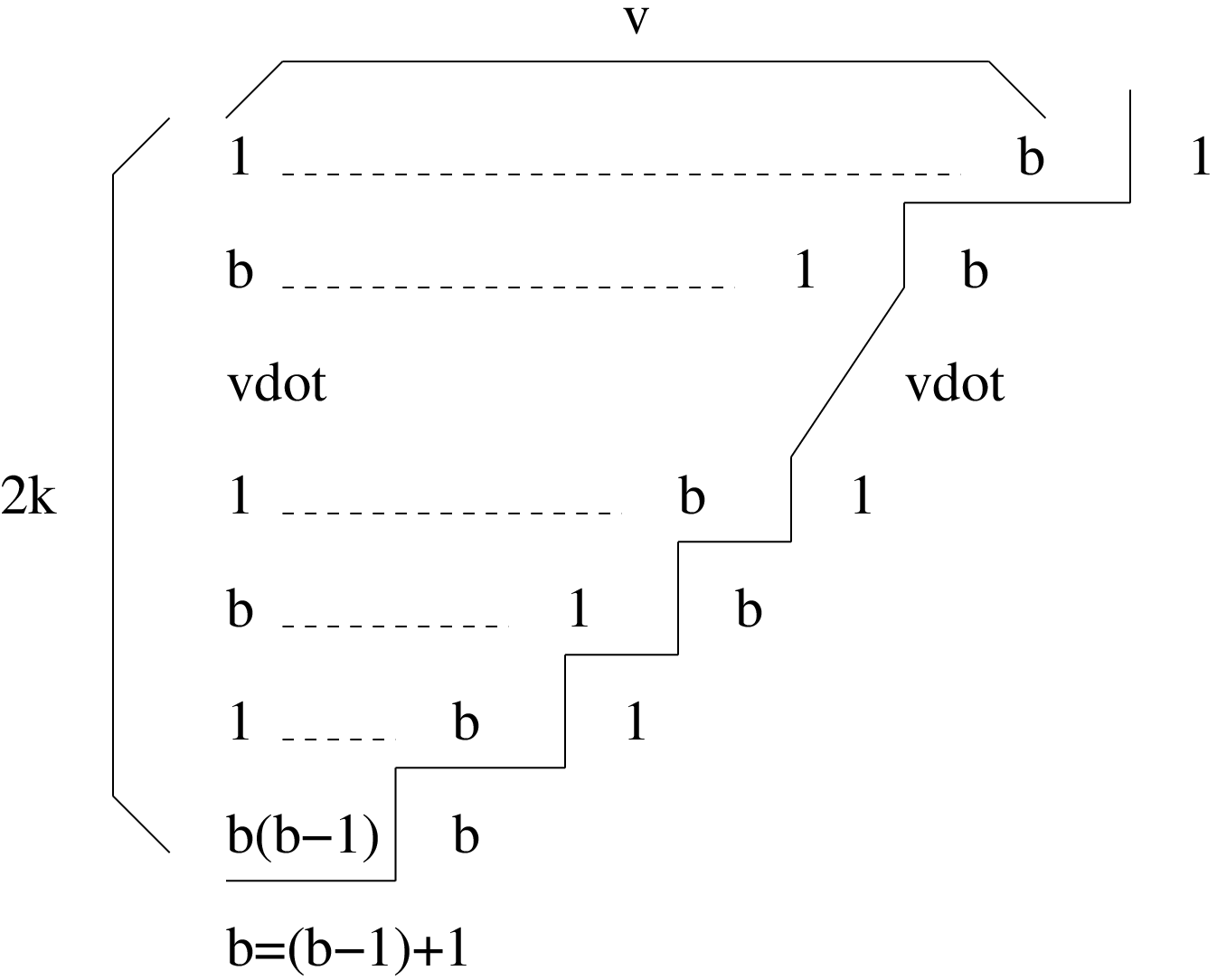}
\caption{Case (6)}\label{figcas6}
\end{center}
\end{figure}

\item [Case (6)] If  $p=(1b)^kb$, let us consider $v=\Phi^{-1}((1b)^k\cdot (b-1))$. This word is illustrated at the left-side of the line in Figure \ref{figcas6}. Only the first and the last letters of each line are written. Since $1$ and $b$ are odd and $(b-1)$ has even length, the $(2k+1)$ first lines have even length, following that the first and the last letters differ. The first letter at the right-side of Figure \ref{figcas6} are written. Then using Lemma \ref{glemma}, we get 
$$\Phi^{-1}(p)=\Phi^{-1}((1b)^kb)=\Phi^{-1}((1b)^k(b-1))\Phi^{-1}((1b)^{k}1).$$
Since $\Phi^{-1}((1b)^{k}1)$ is prefix of the minimal smooth word, it is smaller than $\Phi^{-1}((1b)^kb)$, we conclude that $\Phi^{-1}(p)$ is not a Lyndon word.
\end{enumerate}
{\flushright \QED}

\begin{lemma} \label{lemmm} The word $\Delta^{-1}_1(m_{\{1,b\}})$ is a smooth Lyndon word over the odd alphabet $\{1< b\}$.
\end{lemma}

\Proof By Proposition \ref{motMin} ii), we know that the minimal smooth word $\Phi^{-1}((1b)^\omega)$ is a Lyndon word. It is then sufficient to prove that this implies that $\Delta^{-1}_1(m_{\{1,b\}})$ is also a Lyndon word. Let us consider a factor $f$ of $m_{\{1,b\}}$. Since $m_{\{1,b\}}$ is a Lyndon word, $m_{\{1,b\}}< f$.  Then, let us compute $F=\Delta^{-1}(f)$ in $\Delta^{-1}_1(m_{\{1,b\}})$. There are $2$ cases to consider. 
\begin{itemize}
\item [i)]  If $m_{\{1,b\}}=11y$ and $f=1by'$, with $y \in \A^\omega, y' \in \A^*$: then from Lemma \ref{lemPos}, $f[1]$ is in an odd position, implying that $f$ starts in an even position. Thus, $F=\Delta^{-1}_1(f)=1b^b\cdots$ and $\Delta^{-1}_1(m_{\{1,b\}})=1b1\cdots$. Consequently, $\Delta^{-1}_1(m_{\{1,b\}})< F$. 

\item [ii)] If $m_{\{1,b\}}=1^bx1y$ and $s=1^bxby'$, with $x,y \in \A^*$, $y \in \A^\omega$: then  from Lemma \ref{blocb}, we know that the $b$ in $s$ is at the begining of a block, so it is at an odd position. Furthermore, by Lemma \ref{lemPos}, the block $1^b$ starts in an even position. Thus, $1^bx$ has odd length. Applying Lemma \ref{lemdelta}, we conclude that $\Delta^{-1}_1(m_{\{1,b\}}) < F$, and consequently, $\Delta^{-1}_1(m_{\{1,b\}})$ is an infinite smooth Lyndon word.\\
\vspace{-1cm}
{\flushright \QED}
\end{itemize}

\begin{proposition} \label{propOuf2} Over the alphabet $\{1<b\}$, with $b$ odd, the only infinite smooth Lyndon word $w$ such that $\Phi(w)$ starts by $11$ is $\Delta^{-1}_1(m_{\{1,b\}})$.
\end{proposition}

\begin{figure}[h]
\begin{center}
\includegraphics[width=4.5cm]{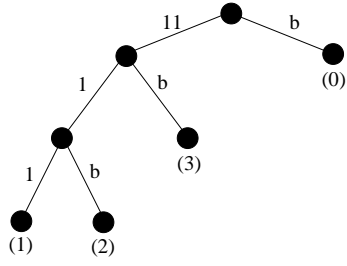}
\caption{Possible prefixes of $p$ starting by $11$}\label{figNon11}
\end{center}
\end{figure}

\Proof Let $p$ be a prefix of $\Phi(w)$. By Lemma \ref{15},  we have $p[0]=1$. Let us prove that the only smooth Lyndon word $w$ such that $p$ is starting by $11$ is $\Delta^{-1}_1(m_{\{1,b\}})$. We proceed by the inspection of the different cases illustrated in Figure \ref{figNon11}.

\begin{itemize}
\item [Case (1)] If $p=1111$, then \\
\indent $\displaystyle \Delta^0(w)=1(b^b\underline 1^b)^\frac{b-1}{2}b^b1b \cdots \\
\indent \Delta^1(w)=1b^b1\cdots\\
\indent \Delta^2(w)=1b \cdots\\
\indent \Delta^3(w)= 1 \cdots$\\
 $w$ has the factor $1^b$.

\item [Case (2)] If $p=111b$, then \\
\indent $\displaystyle \Delta^0(w)=(1b^b)^\frac{b-1}{2}1(b^b\underline 1^b)^\frac{b-1}{2}b^b 1 \cdots \\
\indent \Delta^1(w)=(1b)^\frac{b-1}{2}1b^b1 \cdots\\
\indent \Delta^2(w)=1^bb \cdots\\
\indent \Delta^3(w)=b \cdots$\\
$w$ has the factor $1^b$. 

\item [Case (3)] If $p=11b$, then let us prove that the only smooth Lyndon word $w$ having $p$ has prefix of $\Phi(w)$ is $\Delta^{-1}_1(m_{\{1,b\}})$. To do so, we will prove that $\Phi^{-1}(11bs)$ is a Lyndon word if and only if $\Phi^{-1}(1bs)$ is so.  Let us consider $w=\Phi^{-1}(1bs)$ and $f$, a factor of it. Let us suppose that $w$ is not a Lyndon word and let us show that it implies that $\Delta^{-1}_1(w)$ is not too. Let us rewrite $w=1^bxb y$ and $f=1^bx1 y'$, with $x,y' \in \A^*$ and$y\in \A^\omega$. There are $2$ cases to consider. 
\begin{itemize}
\item [i)] If $x$ ends by $b$: rewrite $f=1^bx'b1y$ and $w=1^bx'bby'$. Then, $w$ contains the factor $b^b$. The existence of this factor in $w$ implies that there is a factor $1^b$ in the word $\Delta^{-1}_1(w)$. Consequently, since $\Delta^{-1}_1(w)$ starts by $1b$, it can not be a Lyndon word.

\item [ii)] If $x$ ends by $1$: rewrite $f=1^bx'11y$ and $w=1^bx'1by'$. By Lemma \ref{lemPos}, $f$ starts in an even position. Then let us compute $\Delta^{-1}_1(f)$ and $\Delta^{-1}_1(w)$. Again, by Lemma \ref{lemPos}, we have that the $b$ in $w$ is in an odd position, implying that $1^bx'1$ has odd length. Using Lemma \ref{lemdelta}, we conclude that $f< w  \Rightarrow \Delta^{-1}_1(f)< \Delta^{-1}_1(w)$: $\Delta^{-1}_1(w)$ is not a Lyndon word.  
\end{itemize}
We conclude, using Lemma \ref{lemmm}.\\
\vspace{-1cm}
{\flushright \QED}
\end{itemize}

\begin{theorem}\label{thmOuf} Over the alphabet $\{1<b\}$, with $b$ odd, the only infinite smooth words that are also  Lyndon words are  $m_{\{1<b\}}$ and  $\Delta^{-1}_1{m_{\{1<b\}}}$.
\end{theorem}

\Proof Follows from the union of Propositions \ref{propOuf1} and \ref{propOuf2}. \QED

\Proof (of Theorem \ref{thmii}) Follows from Theorems \ref{thmii1} and \ref{thmOuf}.\QED

\subsection{Over $\A$ with $a$ odd and  $b$ even}

In this section, we consider infinite smooth words over an alphabet $\{a<b\}$, with $a$ odd and $b$ even.  We prove that over this alphabet, there is no infinite smooth word that is also a Lyndon word. In order to prove it, we consider $2$ cases, $a\neq 1$ and $a=1$, that have to be analysed separately.

\begin{theorem}\label{celui} Over the alphabet $\{a<b\}$, with $a\neq 1$ odd and $b$ even, there is no smooth infinite word that is a Lyndon word.
\end{theorem}

\begin{figure}[ht]
\begin{center}
  \subfigure[for an odd-even alphabet $\{a<b\}$, with $a\neq 1$.]{
    \label{figip}
    \includegraphics[width=4.5cm,keepaspectratio=true]{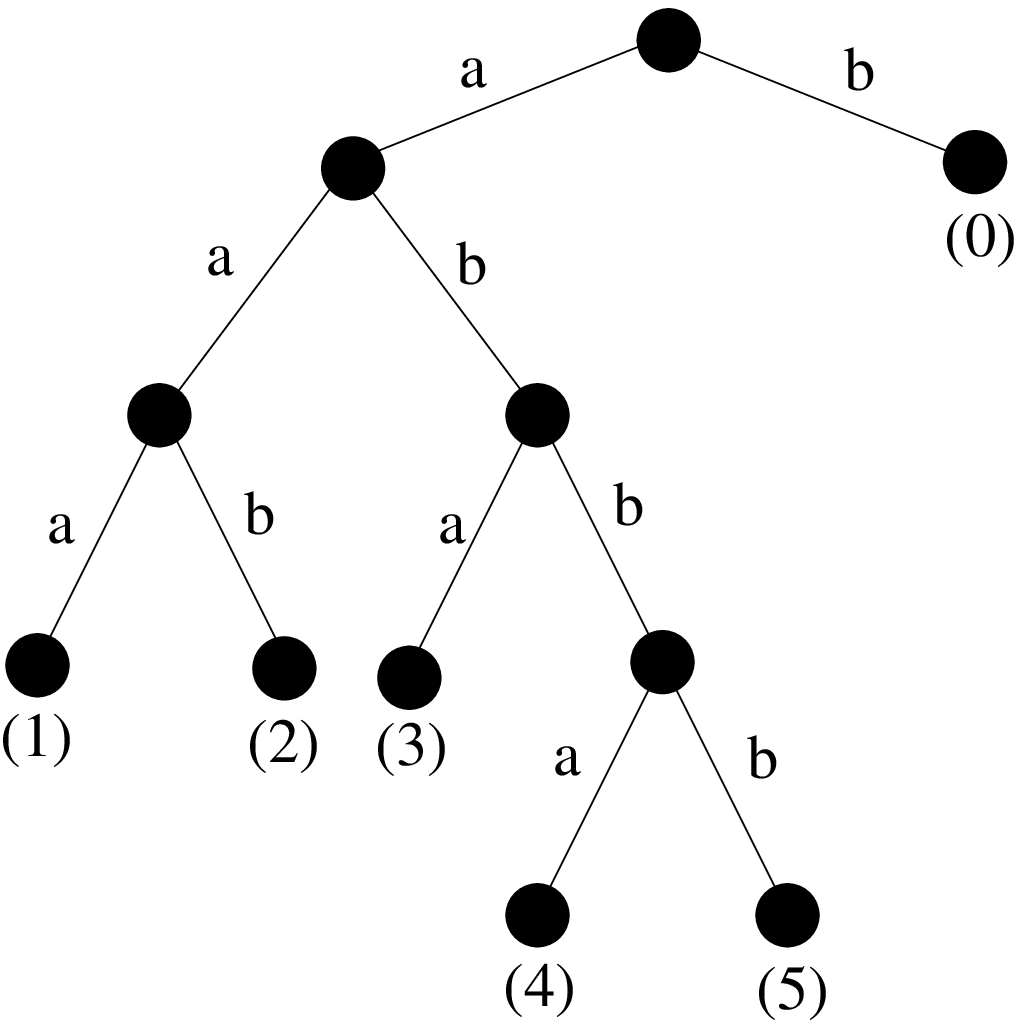}}
  \hspace{1cm}
  \subfigure[for the alphabet $\{a<b\}$, with $a=1$ and $b=4n$.]{
    \label{fig14n}
    \includegraphics[width=5.5cm,keepaspectratio=true]{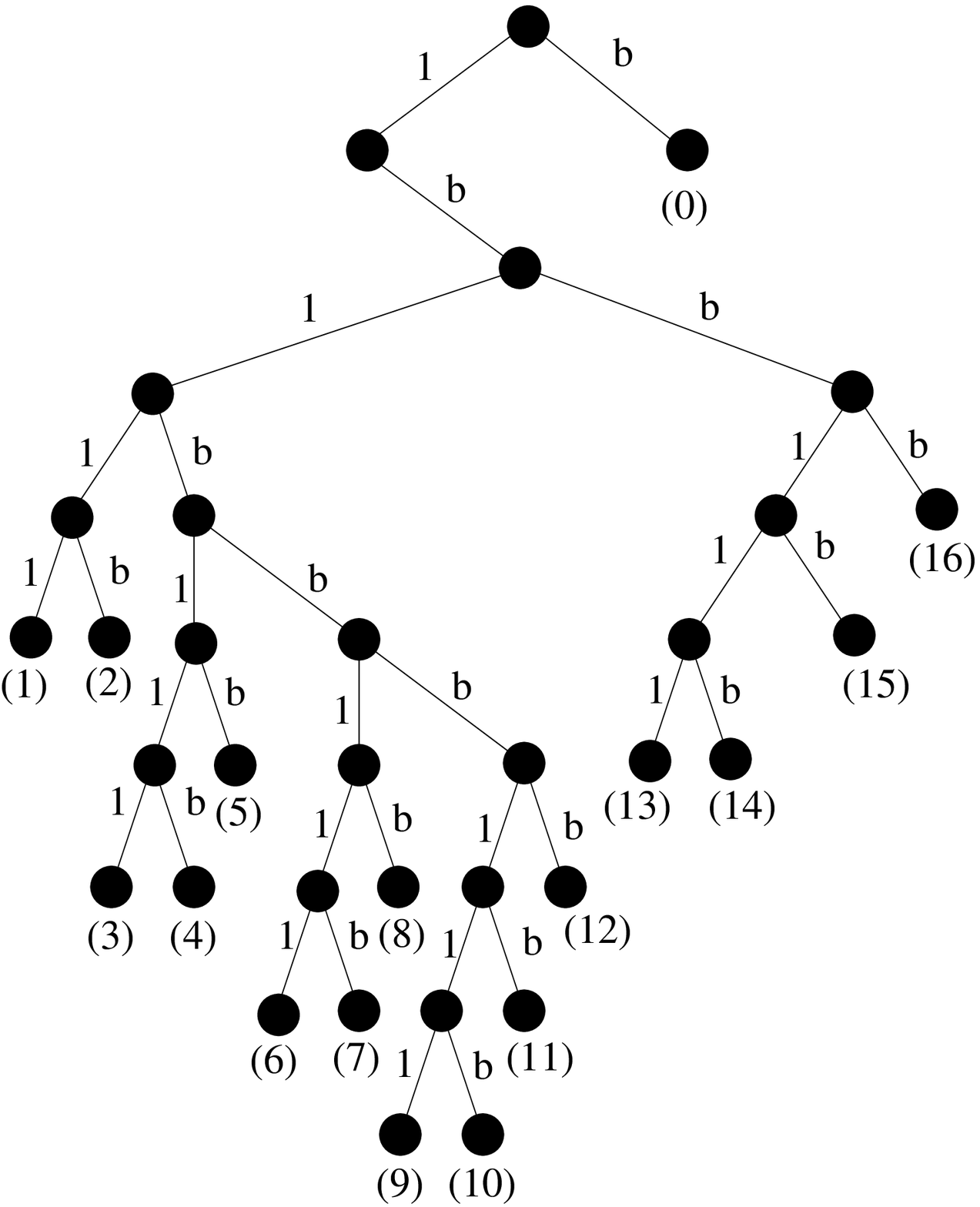}}
  \caption{ \label{doublefig}Possible cases}
\end{center}
\end{figure}

\Proof There are 5 possibilities to consider, illustrated in Figure \ref{doublefig} (a).

\begin{enumerate}
\item [Case (1)] If $p=aaa$, then,\\
\indent $\displaystyle \Delta^0(w)=(a^ab^a)^{\frac{a-1}{2}}a^a(b^b\underline a^b)^\frac{a-1}{2}b^b\cdots \\
\indent \Delta^1(w)=a^ab^a \cdots \\
\indent \Delta^2(w)=aa \cdots$\\
$w$ has the factor $a^b$.

\item [Case (2)] If $p=aab$, then,\\
\indent $\displaystyle \Delta^0(w)=(a^ab^a)^\frac{b}{2}(\underline a^bb^b)^\frac{b}{2}\cdots\\
\indent \Delta^1(w)=a^bb^b\cdots\\
\indent  \Delta^2(w)=bb \cdots $\\
$w$ has the factor $a^b$.

\item [Case (3)] If $p=aba$, then, \\
\indent $\displaystyle \Delta^0(w)=(a^bb^b)^\frac{a-1}{2}\underline a^b(b^aa^a)^\frac{a-1}{2}b^a\cdots \\
\indent \Delta^1(w)=b^a a^a\cdots \\
\indent \Delta^2(w)=aa\cdots$\\
$w$ has  the factor $a^bb^a a^a$.

\item	[Case (4)] If $p=abba$, then, \\
\indent $\displaystyle \Delta^0(w)=((a^bb^b)^\frac{b}{2}(a^ab^a)^\frac{b}{2})^\frac{a-1}{2}(a^bb^b)^\frac{b}{2}((a^ab^a)^\frac{a-1}{2}a^a(b^ba^b)^\frac{a-1}{2}b^b)^\frac{a-1}{2}(a^ab^a)^\frac{a-1}{2}a^a\\
\indent \qquad \qquad  ((b^b\underline a^b)^\frac{b}{2}(b^aa^a)^\frac{b}{2})^\frac{a-1}{2}(b^ba^b)^{\frac{b}{2}}\cdots \\
\indent \Delta^1(w)=(b^ba^b)^\frac{a-1}{2}b^b(a^ab^a)^\frac{a-1}{2}a^a(b^ba^b)^\frac{a-1}{2}b^b \cdots\\
\indent \Delta^2(w)=b^aa^ab^a\cdots \\
\indent  \Delta^3(w)=aa a\cdots$\\
$w$ has the factor $a^b b^aa^a$.

\item	[Case (5)] If $p=abbb$, then, \\
\indent $\displaystyle \Delta^0(w)=((a^bb^b)^\frac{b}{2}(a^ab^a)^\frac{b}{2})^\frac{b}{2}((a^bb^b)^\frac{a-1}{2}\underline a^b(b^aa^a)^\frac{a-1}{2}b^a)^\frac{b}{2}\cdots \\
\indent \Delta^1(w)=(b^ba^b)^\frac{b}{2}(b^a a^a)^\frac{b}{2}\cdots \\
\indent \Delta^2(w)=b^ba^b\cdots\\
\indent \Delta^3(w)=bb\cdots$\\
$w$ has the factor $a^bb^a a^a$.
\end{enumerate}

In each case, it is possible to find a factor $f$ smaller than the smooth word $w$. Thus, there is no smooth Lyndon  word. \QED

\begin{lemma} \label{lem1b} Let $w$ be an infinite smooth Lyndon word over the alphabet $\{1<b\}$, with $b$ even. Let $p$ be a prefix of $\Phi(w)$. Then $p[0,1]=1b$. 
\end{lemma}

\begin{figure}[h]
\begin{center}
\includegraphics[width=3.5cm]{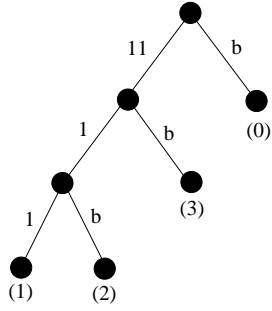}
\caption{Possible cases. }\label{fig99}
\end{center}
\end{figure}

\Proof There are $4$ cases to consider, illustrated in Figure \ref{fig99}. We know from Lemma \ref{15} that the Case (0) is excluded. For the $3$ other cases, we show that $1^b$ is factor of $w$, so it is not a smooth Lyndon word.

\begin{enumerate}
\item [Case (1)] If $p=1111$, then,\\
\indent $\displaystyle \Delta^0(w)=1(b^b\underline 1^b)^\frac{b}{2}b1\cdots \\
\indent \Delta^1(w)=1b^b1\cdots \\
\indent \Delta^2(w)=1b\cdots \\
\indent \Delta^3(w)=1 \cdots$

\item [Case (2)] If $p=111b$, then,\\
\indent $\displaystyle \Delta^0(w)=(1b^b)^\frac{b}{2}(1b)^\frac{b}{2}\underline 1^bb\cdots\\
\indent \Delta^1(w)=(1b)^\frac{b}{2}1^bb\cdots\\
\indent \Delta^2(w)=1^bb\cdots\\
\indent  \Delta^3(w)=b \cdots $

\item [Case (3)] If $p=11b$, then, \\
\indent $\displaystyle \Delta^0(w)=(1b)^\frac{b}{2}\underline 1^bb\cdots \\
\indent \Delta^1(w)=1^bb\cdots \\
\indent \Delta^2(w)=b\cdots$
\end{enumerate}
\vspace{-1cm}
{\flushright \QED}

\begin{theorem}\label{thm14n} Over the alphabet $\{a<b\}$, with $a=1$ and $b=4n$, there is no infinite smooth word that is a Lyndon word. 
\end{theorem}

\Proof Figure \ref{doublefig} (b) shows the different cases to consider. Lemma \ref{lem1b} is used to eliminate the case of a prefix $p$ starting by $11$. For each of the $16$ cases, it is again possible to find a factor of $\Phi^{-1}(p)$ in order to prove that it is not a prefix of an infinite Lyndon word. The details are given in Appendix A. 
\vspace{-0.5cm}
{\flushright \QED}

\begin{theorem}\label{thm122n} Over the alphabet $\{a<b\}$, with $a=1$ and $b=2(2n+1)$, there is no infinite smooth word that is a Lyndon word. 
\end{theorem}

\begin{figure}[h]
\begin{center}
\includegraphics[width=9cm]{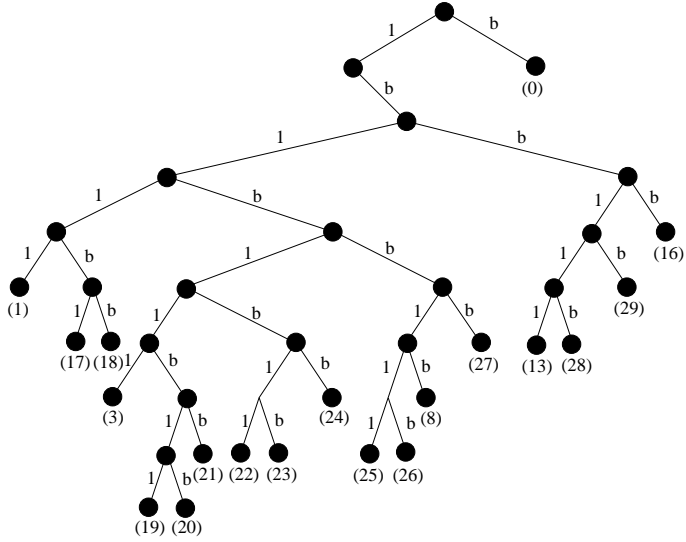}
\caption{Different cases for the alphabet $\{a<b\}$, with $a=1$ and $b=2(2n+1)$}\label{fig122n}
\end{center}
\end{figure}

\Proof Figure \ref{fig122n} shows the different cases to consider.  Cases numbered less or equal to $16$ are the same as in Theorem \ref{thm14n}. For the other cases, it is possible to find a factor in $\Phi^{-1}(p)$ smaller than its prefix, implying that the word is not a Lyndon word. The details are given in Appendix B. Again, we use Lemma \ref{lem1b} to eliminate the case of a prefix $p$ starting by $11$.
\vspace{-0.8cm}
{\flushright \QED}

\section{Summary and concluding remarks}

The next theorem summarizes the results of Section 3.

\begin{theorem} Over any $2$-letter alphabet, the only infinite smooth words that are also infinite Lyndon words are $m_{\{2a<2b\}}$, $m_{\{1<2b+1\}}$ and $\Delta^{-1}(m_{\{1<2b+1\}})$, for $a,b \in \N \setminus \{0\}$. 
\end{theorem}

Recall that for the alphabet $\{1,2\}$, it is conjectured in \cite{bdlv2006} that any smooth factor appears in any infinite smooth word. A direct corollary of this conjecture is that no infinite smooth word is a Lyndon word: there exists a minimal suffix (see \cite{bpm2007}), even smaller than the minimal smooth word and a prefix of this suffix must appear in any smooth word. Since we proved in this paper that there is no smooth Lyndon word over the alphabet $\{1,2\}$, our result  reinforce the conjecture from \cite{bdlv2006}. 

The existence of infinite smooth Lyndon words over the alphabets $\{2a<2b\}$ and $\{1<2b+1\}$ leads to the following corollary.

\begin{corollary} Let $\A$ be a $2$-letter alphabet such that $\A = \{2a<2b\}$ or $\A=\{1<2b+1\}$. Then, any infinite smooth words $w \in \A^\omega$ does not contain every smooth factors. 
\end{corollary}

It is also interesting to notice that our main result completely characterized the trivial finite Lyndon factorization of infinite smooth words: the only infinite smooth words that have a finite Lyndon factorization composed of only one factor are $m_{\{2a<2b\}}$, $m_{\{1<2b+1\}}$ and $\Delta^{-1}(m_{\{1<2b+1\}})$. It is still an open problem to characterized infinite smooth words that have a non trivial finite Lyndon factorization. Giving  an explicit computation of the Lyndon factorization, finite or infinite, of any infinite smooth words, as Melan\c con did for characteristic Sturmian words \cite{gm2000}, is still a challenging problem.

Finally, it would be really nice to find a simpler proof of our result. \\

{\noindent\bf Acknowledgements}. 
The author would like to thank Pierre Lalonde for his interesting question during the LaCIM seminar that leads to this paper and Srecko Brlek for his comments.

\bibliographystyle{alpha} 
{\footnotesize
\bibliography{biblio}
}

\newpage

\appendix

\section{Details of the proof of Theorem \ref{thm14n} }

\begin{itemize}
\item [Case (1)] If $p=1b111$, then\\
\indent $\displaystyle \Delta^0(w)= 1^b((b1)^\frac{b}{2}(b^b\underline{ 1}^b)^\frac{b}{2})^\frac{b}{2}b1^bb\cdots \\
\indent \Delta^1(w)= b(1^bb^b)^\frac{b}{2}1b\cdots \\
\indent \Delta^2(w)= 1b^b1\cdots \\
\indent \Delta^3(w)= 1b\cdots \\
\indent \Delta^4(w)= 1\cdots $\\
$w$ has the factor $1^bb1^b$.

\item [Case (2)] If $p=1b11b$, then\\
\indent $\displaystyle \Delta^0(w)= (1^b(b1)^\frac{b}{2}b^b(1b)^\frac{b}{2})^\frac{b}{4}(\underline 1^bb)^\frac{b}{2}(1^bb^b)^\frac{b}{2}1b \cdots \\
\indent \Delta^1(w)=(b1^b)^\frac{b}{2}(b1)^\frac{b}{2}b^b1 \cdots \\
\indent \Delta^2(w)= (1b)^\frac{b}{2}1^bb\cdots \\
\indent \Delta^3(w)= 1^bb\cdots \\
\indent \Delta^4(w)=b \cdots $\\
$w$ has the factor $(1^bb)^\frac{b}{2}1^b$.

\item [Case (3)] If $p=1b1b111$, then\\
\indent $\displaystyle \Delta^0(w)= (1^bb)^\frac{b}{2}((1^bb^b)^\frac{b}{2}1(b^b1^b)^\frac{b}{2}b)^\frac{b}{4}(((1^bb^b)^\frac{b}{2}(1b)^\frac{b}{2})^\frac{b}{2}(\underline 1^bb)^\frac{b}{2})^\frac{b}{2})^\frac{b}{2}(1^bb^b)^\frac{b}{2}(1b^b)^\frac{b}{2}(1b)^\frac{b}{2}1^bb \cdots \\
\indent \Delta^1(w)= (b1)^\frac{b}{2}((b^b1)^\frac{b}{2}((b^b1^b)^\frac{b}{2}(b1)^\frac{b}{2})^\frac{b}{2})^\frac{b}{2}b^b(1b)^\frac{b}{2}1^bb \cdots \\
\indent \Delta^2(w)= 1^b((b1)^\frac{b}{2}(b^b1^b)^\frac{b}{2})^\frac{b}{2}b1^bb\cdots \\
\indent \Delta^3(w)= b(1^bb^b)^\frac{b}{2}1b\cdots \\
\indent \Delta^4(w)= 1b^b1\cdots \\
\indent \Delta^5(w)= 1b\cdots \\
\indent \Delta^6(w)=1 \cdots $\\
$w$ has the factor $(1^bb)^\frac{b}{2}(1^bb^b)^\frac{b}{2}(1b)^\frac{b}{2}$.

\item [Case (4)] If $p=1b1b11b$, then \\
\indent $\displaystyle \Delta^0(w)=((1^bb)^\frac{b}{2}((1^bb^b)^\frac{b}{2}1(b^b1^b)^\frac{b}{2}b)^\frac{b}{4}((1^bb^b)^\frac{b}{2}(1b)^\frac{b}{2})^\frac{b}{2}(1^b(b1)^\frac{b}{2}b^b(1b)^\frac{b}{2})^\frac{b}{4})^\frac{b}{4}\\
\indent \qquad \qquad ((1^bb)^\frac{b}{2}(1^bb^b)^\frac{b}{2}(1b^b)^\frac{b}{2}(1b)^\frac{b}{2})^\frac{b}{4}((\underline 1^bb)^\frac{b}{2}((1^bb^b)^\frac{b}{2}(1b)^\frac{b}{2})^\frac{b}{2})^\frac{b}{2}1^b(b1)^\frac{b}{2}b^b1 \cdots \\
\indent \Delta^1(w)= ((b1)^\frac{b}{2}(b^b1)^\frac{b}{2}(b^b1^b)^\frac{b}{2}(b1^b)^\frac{b}{2})^\frac{b}{4}((b1)^\frac{b}{2}b^b(1b)^\frac{b}{2}1^b)^\frac{b}{4}((b1)^\frac{b}{2}(b^b1^b)^\frac{b}{2})^\frac{b}{2}b1^bb \cdots \\
\indent \Delta^2(w)= (1^b(b1)^\frac{b}{2}b^b(1b)^\frac{b}{2})^\frac{b}{4}(1^bb)^\frac{b}{2}(1^bb^b)^\frac{b}{2}1b\cdots \\
\indent \Delta^3(w)= (b1^b)^\frac{b}{2}(b1)^\frac{b}{2}b^b1\cdots \\
\indent \Delta^4(w)= (1b)^\frac{b}{2}1^bb\cdots \\
\indent \Delta^5(w)= 1^bb\cdots \\
\indent \Delta^6(w)=b \cdots $\\
$w$ has the factor $(1^bb)^\frac{b}{2}((1^bb^b)^\frac{b}{2}(1b)^\frac{b}{2})^\frac{b}{2}$.

\item [Case (5)] If $p=1b1b1b$, then\\
\indent $\displaystyle \Delta^0(w)=((1^bb)^\frac{b}{2}(1^bb^b)^\frac{b}{2}(1b^b)^\frac{b}{2}(1b)^\frac{b}{2})^\frac{b}{4}((\underline 1^bb)^\frac{b}{2}(1^bb^b)^\frac{b}{2}(1b)^\frac{b}{2})^\frac{b}{2}1^b(b1)^\frac{b}{2}b^b1 \cdots \\
\indent \Delta^1(w)= ((b1)^\frac{b}{2}b^b(1b)^\frac{b}{2}1^b)^\frac{b}{4}((b1)^\frac{b}{2}(b^b1^b)^\frac{b}{2})^\frac{b}{2}b1^bb\cdots \\
\indent \Delta^2(w)=(1^bb)^\frac{b}{2}(1^bb^b)^\frac{b}{2}1b \cdots \\
\indent \Delta^3(w)=(b1)^\frac{b}{2}b^b1 \cdots \\
\indent \Delta^4(w)=1^bb \cdots \\
\indent \Delta^5(w)=b \cdots $\\
$w$ has the factor $(1^bb)^\frac{b}{2}(1^bb^b)^\frac{b}{2}(1b)^\frac{b}{2}$.

\item [Case (6)] If $p=1b1bb111$, then\\
\indent $\displaystyle \Delta^0(w)= ((1^bb)^\frac{b}{2}((1^bb^b)^\frac{b}{2}(1b)^\frac{b}{2})^\frac{b}{2})^\frac{b}{2}((1^b((b1)^\frac{b}{2}(b^b1^b)^\frac{b}{2})^\frac{b}{2}b((1^bb^b)^\frac{b}{2}(1b)^\frac{b}{2})^\frac{b}{2})^\frac{b}{4} \\
\indent \qquad \qquad ((1^b(b1)^\frac{b}{2}b^b(1b)^\frac{b}{2})^\frac{b}{4}((1^bb)^\frac{b}{2}((1^bb^b)^\frac{b}{2}(1b)^\frac{b}{2})^\frac{b}{2})^\frac{b}{2})^\frac{b}{2})^\frac{b}{2}1^b(((b1)^\frac{b}{2}(b^b\underline1^b)^\frac{b}{2})^\frac{b}{2}(b1^b)^\frac{b}{2})^\frac{b}{2}(b1)^\frac{b}{2}b^b1 \cdots \\
\indent \Delta^1(w)= ((b1)^\frac{b}{2}(b^b1^b)^\frac{b}{2})^\frac{b}{2}((b(1^bb^b)^\frac{b}{2}1(b^b1^b)^\frac{b}{2})^\frac{b}{4}((b1^b)^\frac{b}{2}((b1)^\frac{b}{2}(b^b1^b)^\frac{b}{2})^\frac{b}{2})^\frac{b}{2})^\frac{b}{2}b((1^bb^b)^\frac{b}{2}(1b)^\frac{b}{2})^\frac{b}{2}1^bb \cdots \\
\indent \Delta^2(w)= (1^bb^b)^\frac{b}{2}((1b^b)^\frac{b}{2}((1b)^\frac{b}{2}(1^bb^b)^\frac{b}{2})^\frac{b}{2})^\frac{b}{2}1(b^b1^b)^\frac{b}{2}b1\cdots \\
\indent \Delta^3(w)= b^b((1b)^\frac{b}{2}(1^bb^b)^\frac{b}{2})^\frac{b}{2}1b^b1\cdots \\
\indent \Delta^4(w)= b(1^bb^b)^\frac{b}{2}1b\cdots \\
\indent \Delta^5(w)= 1b^b1\cdots \\
\indent \Delta^6(w)= 1b\cdots \\
\indent \Delta^7(w)= 1\cdots $\\
$w$ has the factor $1^b(b1^b)^\frac{b}{2}(b1)^\frac{b}{2}$ in $(((b1)^\frac{b}{2}(b^b\underline1^b)^\frac{b}{2})^\frac{b}{2}(b1^b)^\frac{b}{2})^\frac{b}{2}$.

\item [Case (7)] If  $p=1b1bb11b$, then\\
\indent $\displaystyle \Delta^0(w)= (((1^bb)^\frac{b}{2}((1^bb^b)^\frac{b}{2}(1b)^\frac{b}{2})^\frac{b}{2})^\frac{b}{2}(1^b((b1)^\frac{b}{2}(b^b1^b)^\frac{b}{2})^\frac{b}{2}b((1^bb^b)^\frac{b}{2}(1b)^\frac{b}{2})^\frac{b}{2})^\frac{b}{4}\\
\indent \qquad \qquad (1^b(b1)^\frac{b}{2}b^b(1b)^\frac{b}{2})^\frac{b}{4}((1^bb)^\frac{b}{2}(1^bb^b)^\frac{b}{2}(1b^b)^\frac{b}{2}(1b)^\frac{b}{2})^\frac{b}{2})^\frac{b}{4}\\
\indent \qquad \qquad  ((1^bb)^\frac{b}{2}((1^bb^b)^\frac{b}{2}(1b)^\frac{b}{2})^\frac{b}{2})^\frac{b}{2}1^b(((b1)^\frac{b}{2}(b^b\underline 1^b)^\frac{b}{2})^\frac{b}{2}(b1^b)^\frac{b}{2})^\frac{b}{2}(b1)^\frac{b}{2}\cdots \\
\indent \Delta^1(w)= (((b1)^\frac{b}{2}(b^b1^b)^\frac{b}{2})^\frac{b}{2}(b(1^bb^b)^\frac{b}{2}1(b^b1^b)^\frac{b}{2})^\frac{b}{4}(b1^b)^\frac{b}{2}((b1)^\frac{b}{2}b^b(1b)^\frac{b}{2}1^b)^\frac{b}{4})^\frac{b}{4}\\
\indent \qquad \qquad  (((b1)^\frac{b}{2}(b^b1^b)^\frac{b}{2})^\frac{b}{2}b((1^bb^b)^\frac{b}{2}(1b)^\frac{b}{2})^\frac{b}{2}1^b)^\frac{b}{4}(((b1)^\frac{b}{2}(b^b1^b)^\frac{b}{2})^\frac{b}{2}(b1^b)^\frac{b}{2})^\frac{b}{2}(b1)^\frac{b}{2}b^b1b \cdots \\
\indent \Delta^2(w)= ((1^bb^b)^\frac{b}{2}(1b^b)^\frac{b}{2}(1b)^\frac{b}{2}(1^bb)^\frac{b}{2})^\frac{b}{4} ((1^bb^b)^\frac{b}{2}1(b^b1^b)^\frac{b}{2}b)^\frac{b}{4}((1^bb^b)^\frac{b}{2}(1b)^\frac{b}{2})^\frac{b}{2}1^bb1 \cdots \\
\indent \Delta^3(w)= (b^b(1b)^\frac{b}{2}1^b(b1)^\frac{b}{2})^\frac{b}{4}(b^b1)^\frac{b}{2}(b^b1^b)^\frac{b}{2}b1\cdots \\
\indent \Delta^4(w)= (b1^b)^\frac{b}{2}(b1)^\frac{b}{2}b^b1\cdots \\
\indent \Delta^5(w)= (1b)^\frac{b}{2}1^bb\cdots \\
\indent \Delta^6(w)= 1^bb\cdots \\
\indent \Delta^7(w)=b \cdots $\\
$w$ has the factor $1^b(b1^b)^\frac{b}{2}(b1)^\frac{b}{2}$. 

\item [Case (8)] If $p=1b1bb1b$, then\\
\indent $\displaystyle \Delta^0(w)= (((1^bb)^\frac{b}{2}((1^bb^b)^\frac{b}{2}(1b)^\frac{b}{2})^\frac{b}{2})^\frac{b}{2}1^b(((b1)^\frac{b}{2}(b^b\underline 1^b)^\frac{b}{2})^\frac{b}{2}(b1^b)^\frac{b}{2})^\frac{b}{2}(b1)^\frac{b}{2}\cdots \\
\indent \Delta^1(w)= (((b1)^\frac{b}{2}(b^b1^b)^\frac{b}{2})^\frac{b}{2}b((1^bb^b)^\frac{b}{2}(1b)^\frac{b}{2})^\frac{b}{2}1^b)^\frac{b}{4}(((b1)^\frac{b}{2}(b^b1^b)^\frac{b}{2})^\frac{b}{2}(b1^b)^\frac{b}{2})^\frac{b}{2}(b1)^\frac{b}{2}b^b1\cdots \\
\indent \Delta^2(w)= ((1^bb^b)^\frac{b}{2}1(b^b 1^b)^\frac{b}{2}b)^\frac{b}{4}((1^bb^b)^\frac{b}{2}(1b)^\frac{b}{2})^\frac{b}{2}1^bb \cdots\\
\indent \Delta^3(w)=(b^b1)^\frac{b}{2}(b^b1^b)^\frac{b}{2}b1 \cdots \\
\indent \Delta^4(w)=(b1)^\frac{b}{2}b^b1 \cdots \\
\indent \Delta^5(w)= 1^bb\cdots \\
\indent \Delta^6(w)=b \cdots $\\
$w$ has the factor $1^b(b1^b)^\frac{b}{2}(b1)^\frac{b}{2}$.

\item [Case (9)] If $p=1b1bbb111$, then\\
\indent $\displaystyle \Delta^0(w)= (((1^bb)^\frac{b}{2}((1^bb^b)^\frac{b}{2}(1b)^\frac{b}{2})^\frac{b}{2})^\frac{b}{2}(1^b(b1)^\frac{b}{2}b^b(1b)^\frac{b}{2})^\frac{b}{4})^\frac{b}{2}(((1^bb)^\frac{b}{2}((1^bb^b)^\frac{b}{2}1(b^b1^b)^\frac{b}{2}b)^\frac{b}{4}\\
\indent \qquad \qquad ((1^bb^b)^\frac{b}{2}(1b)^\frac{b}{2})^\frac{b}{2}(1^b(b1)^\frac{b}{2}b^b(1b)^\frac{b}{2})^\frac{b}{4})^\frac{b}{4}(((1^bb)^\frac{b}{2}(1^bb^b)^\frac{b}{2}(1b^b)^\frac{b}{2}(1b)^\frac{b}{2})^\frac{b}{4}\\
\indent \qquad \qquad (((\underline 1^bb)^\frac{b}{2}((1^bb^b)^\frac{b}{2}(1b)^\frac{b}{2})^\frac{b}{2})^\frac{b}{2}(1^b(b1)^\frac{b}{2}b^b(1b)^\frac{b}{2})^\frac{b}{4})^\frac{b}{2})^\frac{b}{2})^\frac{b}{2}(1^bb)^\frac{b}{2}(((1^bb^b)^\frac{b}{2}1(b^b1^b)^\frac{b}{2}b)^\frac{b}{4}\\
\indent \qquad \qquad (((1^bb^b)^\frac{b}{2}(1b)^\frac{b}{2})^\frac{b}{2}(1^bb)^\frac{b}{2})^\frac{b}{2})^\frac{b}{2}(1^bb^b)^\frac{b}{2}(1b^b)^\frac{b}{2}(1b)^\frac{b}{2}1^bb \cdots \\
\indent \Delta^1(w)= (((b1)^\frac{b}{2}(b^b1^b)^\frac{b}{2})^\frac{b}{2}(b1^b)^\frac{b}{2})^\frac{b}{2}(((b1)^\frac{b}{2}(b^b1)^\frac{b}{2}(b^b1^b)^\frac{b}{2}(b1^b)^\frac{b}{2})^\frac{b}{4}(((b1)^\frac{b}{2}b^b(1b)^\frac{b}{2}1^b)^\frac{b}{4}\\
\indent \qquad \qquad (((b1)^\frac{b}{2}(b^b1^b)^\frac{b}{2})^\frac{b}{2}(b1^b)^\frac{b}{2})^\frac{b}{2})^\frac{b}{2})^\frac{b}{2}(b1)^\frac{b}{2}((b^b1)^\frac{b}{2}((b^b1^b)^\frac{b}{2}(b1)^\frac{b}{2})^\frac{b}{2})^\frac{b}{2}b^b(1b)^\frac{b}{2}1^bb \cdots \\
\indent \Delta^2(w)= ((1^bb^b)^\frac{b}{2}(1b)^\frac{b}{2})^\frac{b}{2}((1^b(b1)^\frac{b}{2}b^b(1b)^\frac{b}{2})^\frac{b}{4}((1^bb)^\frac{b}{2}((1^bb^b)^\frac{b}{2}(1b)^\frac{b}{2})^\frac{b}{2})^\frac{b}{2})^\frac{b}{2}1^b((b1)^\frac{b}{2}(b^b1^b)^\frac{b}{2})^\frac{b}{2}b1^bb\cdots \\
\indent \Delta^3(w)= (b^b1^b)^\frac{b}{2}((b1^b)^\frac{b}{2}((b1)^\frac{b}{2}(b^b1^b)^\frac{b}{2})^\frac{b}{2})^\frac{b}{2}b(1^bb^b)^\frac{b}{2}1b \cdots \\
\indent \Delta^4(w)= b^b((1b)^\frac{b}{2}(1^bb^b)^\frac{b}{2})^\frac{b}{2}1b^b1\cdots \\
\indent \Delta^5(w)= b(1^bb^b)^\frac{b}{2}1b\cdots \\
\indent \Delta^6(w)= 1b^b1\cdots \\
\indent \Delta^7(w)= 1b\cdots \\
\indent \Delta^8(w)=1 \cdots $\\
$w$ has the factor\\
 $(((1^bb)^\frac{b}{2}((1^bb^b)^\frac{b}{2}(1b)^\frac{b}{2})^\frac{b}{2})^\frac{b}{2}(1^b(b1)^\frac{b}{2}b^b(1b)^\frac{b}{2})^\frac{b}{4})^\frac{b}{2}(1^bb)^\frac{b}{2}(((1^bb^b)^\frac{b}{2}1(b^b1^b)^\frac{b}{2}b)^\frac{b}{4} (((1^bb^b)^\frac{b}{2}(1b)^\frac{b}{2})^\frac{b}{2}(1^bb)^\frac{b}{2})^\frac{b}{2})^\frac{b}{2}$.

\item [Case (10)] If $p=1b1bbb11b$, then\\
\indent $\displaystyle \Delta^0(w)= (((1^bb)^\frac{b}{2}((1^bb^b)^\frac{b}{2}(1b)^\frac{b}{2})^\frac{b}{2})^\frac{b}{2}(1^b(b1)^\frac{b}{2}b^b(1b)^\frac{b}{2})^\frac{b}{4})^\frac{b}{2}((1^bb)^\frac{b}{2}((1^bb^b)^\frac{b}{2}1(b^b1^b)^\frac{b}{2}b)^\frac{b}{4}\\
\indent \qquad \qquad ((1^bb^b)^\frac{b}{2}(1b)^\frac{b}{2})^\frac{b}{2}(1^b(b1)^\frac{b}{2}b^b(1b)^\frac{b}{2})^\frac{b}{4})^\frac{b}{4}(( 1^bb)^\frac{b}{2}(1^bb^b)^\frac{b}{2}(1b^b)^\frac{b}{2}(1b)^\frac{b}{2})^\frac{b}{4}(((\underline1^bb)^\frac{b}{2}((1^bb^b)^\frac{b}{2}(1b)^\frac{b}{2})^\frac{b}{2})^\frac{b}{2}1^b\\
\indent \qquad\qquad  (((b1)^\frac{b}{2}(b^b1^b)^\frac{b}{2})^\frac{b}{2}(b1^b)^\frac{b}{2}(b1)^\frac{b}{2}\cdots \\
\indent \Delta^1(w)= ((((b1)^\frac{b}{2}(b^b1^b)^\frac{b}{2})^\frac{b}{2}(b1^b)^\frac{b}{2})^\frac{b}{2}((b1)^\frac{b}{2}(b^b1)^\frac{b}{2}(b^b1^b)^\frac{b}{2}(b1^b)^\frac{b}{2})^\frac{b}{4}((b1)^\frac{b}{2}b^b(1b)^\frac{b}{2}1^b)^\frac{b}{4}\\
\indent \qquad \qquad  (((b1)^\frac{b}{2}(b^b1^b)^\frac{b}{2})^\frac{b}{2}b((1^bb^b)^\frac{b}{2}(1b)^\frac{b}{2})^\frac{b}{2}1^b)^\frac{b}{4})^\frac{b}{4}(((b1^b)^\frac{b}{2}((b1)^\frac{b}{2}(b^b1^b)^\frac{b}{2})^\frac{b}{2})^\frac{b}{2}b(((1^bb^b)^\frac{b}{2}(1b)^\frac{b}{2})^\frac{b}{2}(1^bb)^\frac{b}{2})^\frac{b}{2}\\
\indent \qquad \qquad  (1^bb^b)^\frac{b}{2})^\frac{b}{4}(((1b^b)^\frac{b}{2}((1b)^\frac{b}{2}(1^bb^b)^\frac{b}{2})^\frac{b}{2})^\frac{b}{2}(1(b^b1^b)^\frac{b}{2}b(1^bb^b)^\frac{b}{2})^\frac{b}{4})^\frac{b}{2}(1b^b)^\frac{b}{2}(1b)^\frac{b}{2}1^bb1 \cdots \\
\indent \Delta^2(w)= (((1^bb^b)^\frac{b}{2}(1b)^\frac{b}{2})^\frac{b}{2}(1^b(b1)^\frac{b}{2}b^b(1b)^\frac{b}{2})^\frac{b}{4}(1^bb)^\frac{b}{2}((1^bb^b)^\frac{b}{2}1(b^b1^b)^\frac{b}{2}b)^\frac{b}{4})^\frac{b}{4}\\
\indent \qquad \qquad  (((1b)^\frac{b}{2}(1^bb^b)^\frac{b}{2})^\frac{b}{2}1((b^b1^b)^\frac{b}{2}(b1)^\frac{b}{2})^\frac{b}{2}b^b)^\frac{b}{4}(((1b)^\frac{b}{2}(1^bb^b)^\frac{b}{2})^\frac{b}{2}(1b^b)^\frac{b}{2})^\frac{b}{2}(1b)^\frac{b}{2}1^bb1
 \cdots \\
\indent \Delta^3(w)= ((b^b1^b)^\frac{b}{2}(b1^b)^\frac{b}{2}(b1)^\frac{b}{2}(b^b1)^\frac{b}{2})^\frac{b}{4}((1^bb^b)^\frac{b}{2}1(b^b1^b)^\frac{b}{2}b)^\frac{b}{4}((1^bb^b)^\frac{b}{2}(1b)^\frac{b}{2})^\frac{b}{2}1^bb1 \cdots \\
\indent \Delta^4(w)= (b^b(1b)^\frac{b}{2}1^b(b1)^\frac{b}{2})^\frac{b}{4}(b^b1)^\frac{b}{2}(b^b1^b)^\frac{b}{2}b1\cdots  $\\
Notice that we take $\Delta^4(w)$ from Case (7). \\
$w$ has the factor $(((1^bb)^\frac{b}{2}((1^bb^b)^\frac{b}{2}(1b)^\frac{b}{2})^\frac{b}{2})^\frac{b}{2}1^b (((b1)^\frac{b}{2}(b^b1^b)^\frac{b}{2})^\frac{b}{2}$.

\item [Case (11)] If $p=1b1bbb1b$, there are $2$ subcases to consider.
\begin{itemize}
\item [(i)] If $b\neq 4$, there are again $2$ subcases.
\begin{itemize}
\item [(a)] If $b/4$ is odd, then\\
\indent $\displaystyle \Delta^0(w)= (((1^bb)^\frac{b}{2}((1^bb^b)^\frac{b}{2}(1b)^\frac{b}{2})^\frac{b}{2})^\frac{b}{2}(1^b(b1)^\frac{b}{2}b^b(1b)^\frac{b}{2})^\frac{b}{4})^\frac{b}{2}(1^bb)^\frac{b}{2}(((1^bb^b)^\frac{b}{2}1(b^b1^b)^\frac{b}{2}b)^\frac{b}{4}\\
\indent \qquad \qquad  ((( ^bb^b)^\frac{b}{2}(1b)^\frac{b}{2})^\frac{b}{2}(1^bb)^\frac{b}{2})^\frac{b}{2})^\frac{b}{2}(1^bb^b)^\frac{b}{2}((((1b^b)^\frac{b}{2}((1b)^\frac{b}{2}(1^bb^b)^\frac{b}{2})^\frac{b}{2})^\frac{b}{2}(1(b^b1^b)^\frac{b}{2}b\\
\indent \qquad \qquad (1^bb^b)^\frac{b}{2})^\frac{b}{4})^\frac{b}{2} (1b^b)^\frac{b}{2}(((1b)^\frac{b}{2}1^b(b1)^\frac{b}{2}b^b)^\frac{b}{4} (((1b)^\frac{b}{2}(1^bb^b)^\frac{b}{2})^\frac{b}{2}(1b^b)^\frac{b}{2})^\frac{b}{2})^\frac{b}{2}(1b)^\frac{b}{2}\\
\indent \qquad \qquad (((1^bb)^\frac{b}{2}((1^bb^b)^\frac{b}{2}(1b)^\frac{b}{2})^\frac{b}{2})^\frac{b}{2}(1^b(b1)^\frac{b}{2}b^b(1b)^\frac{b}{2})^\frac{b}{4})^\frac{b}{2}(1^bb)^\frac{b}{2}(((1^bb^b)^\frac{b}{2}1(b^b\underline1^b)^\frac{b}{2}b)^\frac{b}{4}\\
\indent \qquad \qquad  (((1^bb^b)^\frac{b}{2}(1b)^\frac{b}{2})^\frac{b}{2}( 1^bb)^\frac{b}{2})^\frac{b}{2})^\frac{b}{2}(1^bb^b)^\frac{b}{2})^\frac{b-4}{8}\cdots\\
\indent \Delta^1(w)= (((b1)^\frac{b}{2}(b^b1^b)^\frac{b}{2})^\frac{b}{2}(b1^b)^\frac{b}{2})^\frac{b}{2}(b1)^\frac{b}{2}((b^b1)^\frac{b}{2}((b^b1^b)^\frac{b}{2}(b1)^\frac{b}{2})^\frac{b}{2})^\frac{b}{2}b^b
((((1b)^\frac{b}{2}(1^bb^b)^\frac{b}{2})^\frac{b}{2}\\
\indent \qquad \qquad (1b^b)^\frac{b}{2})^\frac{b}{2}(1b)^\frac{b}{2}((1^bb)^\frac{b}{2}((1^bb^b)^\frac{b}{2}(1b)^\frac{b}{2})^\frac{b}{2})^\frac{b}{2}1^b(((b1)^\frac{b}{2}(b^b1^b)^\frac{b}{2})^\frac{b}{2}(b1^b)^\frac{b}{2})^\frac{b}{2}(b1)^\frac{b}{2}\\
\indent \qquad \qquad ((b^b1)^\frac{b}{2}((b^b1^b)^\frac{b}{2}(b1)^\frac{b}{2})^\frac{b}{2})^\frac{b}{2}b^b)^\frac{b-4}{8}((((1b)^\frac{b}{2}(1^bb^b)^\frac{b}{2})^\frac{b}{2}(1b^b)^\frac{b}{2})^\frac{b}{2}((1b)^\frac{b}{2}1^b(b1)^\frac{b}{2}b^b)^\frac{b}{4})^\frac{b}{2}\\
\indent \qquad \qquad ((1b)^\frac{b}{2}(1^bb^b)^\frac{b}{2})^\frac{b}{2}1(b^b1^b)^\frac{b}{2}b1  \cdots \\
\indent \Delta^2(w)= (((1^bb^b)^\frac{b}{2}(1b)^\frac{b}{2})^\frac{b}{2}1^b((b1)^\frac{b}{2}(b^b1^b)^\frac{b}{2})^\frac{b}{2}b)^\frac{b}{4}(((1^bb^b)^\frac{b}{2}(1b)^\frac{b}{2})^\frac{b}{2}(1^bb)^\frac{b}{2})^\frac{b}{2}(1^bb^b)^\frac{b}{2}1b^b1\cdots \\
\indent \Delta^3(w)= ((b^b1^b)^\frac{b}{2}b(1^bb^b)^\frac{b}{2}1)^\frac{b}{4}((b^b1^b)^\frac{b}{2}(b1)^\frac{b}{2})^\frac{b}{2}b^b1b\cdots \\
\indent \Delta^4(w)= (b^b1)^\frac{b}{2}(b^b1^b)^\frac{b}{2}b1\cdots \\
\indent \Delta^5(w)= (b1)^\frac{b}{2}b^b1\cdots \\
\indent \Delta^6(w)= 1^bb\cdots \\
\indent \Delta^7(w)=b \cdots $\\
$w$ has the factor $(1^bb)(((1^bb^b)^\frac{b}{2}(1b)^\frac{b}{2})^\frac{b}{2}(1^bb)^\frac{b}{2})^\frac{b}{2}(1^bb^b)^\frac{b}{2}$.

\item [(b)] If $b/4$ is even, then\\
\indent $\displaystyle \Delta^0(w)= ((((1^bb)^\frac{b}{2}((1^bb^b)^\frac{b}{2}(1b)^\frac{b}{2})^\frac{b}{2})^\frac{b}{2}(1^b(b1)^\frac{b}{2}b^b(1b)^\frac{b}{2})^\frac{b}{4})^\frac{b}{2}(1^bb)^\frac{b}{2}(((1^bb^b)^\frac{b}{2}1(b^b1^b)^\frac{b}{2}b)^\frac{b}{4}\\
\indent \qquad \qquad (((1^bb^b)^\frac{b}{2}(1b)^\frac{b}{2})^\frac{b}{2}(1^bb)^\frac{b}{2})^\frac{b}{2})^\frac{b}{2}(1^bb^b)^\frac{b}{2}(((1b^b)^\frac{b}{2}((1b)^\frac{b}{2}(1^bb^b)^\frac{b}{2})^\frac{b}{2})^\frac{b}{2}(1(b^b1^b)^\frac{b}{2}\\
\indent \qquad \qquad b(1^bb^b)^\frac{b}{2})^\frac{b}{4})^\frac{b}{2}(1b^b)^\frac{b}{2}(((1b)^\frac{b}{2}1^b(b1)^\frac{b}{2}b^b)^\frac{b}{4} (((1b)^\frac{b}{2}(1^bb^b)^\frac{b}{2})^\frac{b}{2} (1b^b)^\frac{b}{2})^\frac{b}{2})^\frac{b}{2}(1b)^\frac{b}{2})^\frac{b}{8}\\
\indent \qquad \qquad((((1^bb)^\frac{b}{2}((1^bb^b)^\frac{b}{2}(1b)^\frac{b}{2})^\frac{b}{2})^\frac{b}{2}(1^b(b1)^\frac{b}{2}b^b(1b)^\frac{b}{2})^\frac{b}{4})^\frac{b}{2}((1^bb)^\frac{b}{2}(1^bb^b)^\frac{b}{2}\\
\indent \qquad \qquad (1b^b)^\frac{b}{2}(1b)^\frac{b}{2})^\frac{b}{4})^\frac{b}{2}(( \underline 1^bb)^\frac{b}{2}((1^bb^b)^\frac{b}{2}(1b)^\frac{b}{2})^\frac{b}{2})^\frac{b}{2}1^b((b1)^\frac{b}{2}(b^b1^b)^\frac{b}{2})^\frac{b}{2}b1^bb \cdots\\
\indent \Delta^1(w)= ((((b1)^\frac{b}{2}(b^b1^b)^\frac{b}{2})^\frac{b}{2}(b1^b)^\frac{b}{2})^\frac{b}{2}(b1)^\frac{b}{2}((b^b1)^\frac{b}{2}((b^b1^b)^\frac{b}{2}(b1)^\frac{b}{2})^\frac{b}{2})^\frac{b}{2}b^b\\
\indent \qquad \qquad (((1b)^\frac{b}{2}(1^bb^b)^\frac{b}{2})^\frac{b}{2}(1b^b)^\frac{b}{2})^\frac{b}{2}(1b)^\frac{b}{2}((1^bb)^\frac{b}{2}((1^bb^b)^\frac{b}{2}(1b)^\frac{b}{2})^\frac{b}{2})^\frac{b}{2}1^b)^\frac{b}{8}\\
\indent \qquad \qquad ((((b1)^\frac{b}{2}(b^b1^b)^\frac{b}{2})^\frac{b}{2}(b1^b)^\frac{b}{2})^\frac{b}{2}((b1)^\frac{b}{2}b^b(1b)^\frac{b}{2}1^b)^\frac{b}{4})^\frac{b}{2}((b1)^\frac{b}{2}(b^b1^b)^\frac{b}{2})^\frac{b}{2}b(1^bb^b)^\frac{b}{2}1b \cdots$ \\
Notice that the $\Delta^i(w)$'s are the same as in $(a)$ for $2\leq i \leq 7$.\\
$w$ has the factor $((1^bb)^\frac{b}{2}((1^bb^b)^\frac{b}{2}(1b)^\frac{b}{2})^\frac{b}{2})^\frac{b}{2}1^b((b1)^\frac{b}{2}(b^b1^b)^\frac{b}{2})^\frac{b}{2}$.
\end{itemize}

\item [(ii)] If $b=4$, there are again 3 subcases to consider.

\begin{itemize}
\item [(a)] If $p=1414441411$, then\\
\indent $\displaystyle \Delta^0(w)= (((1^44)^2((1^44^4)^2(14)^2)^2)^21^4(41)^24^4(14)^2)^2(1^44)^2((1^44^4)^21(4^41^4)^24(((1^44^4)^2(14)^2)^2\\
\indent \qquad \qquad (1^44)^2)^2)^2(1^44^4)^2((((14^4)^2((14)^2(1^44^4)^2)^2)^21(4^41^4)^24(1^44^4)^2)^2(14^4)^2(14)^2(1^44)^2\\
\indent \qquad \qquad(1^44^4)^2)^2((14^4)^2((14)^2(1^44^4)^2)^2)^2(1((4^41^4)^2(41)^2)^24^4((14)^2(1^44^4)^2)^2(1(4^41^4)^24(1^44^4)^2\\
\indent \qquad \qquad((14^4)^2((14)^2(1^44^4)^2)^2)^2)^2)^21((((((4^41^4)^2(41)^2)^2(4^41)^2)^2(4^41^4)^24(1^44^4)^21)^2\\
\indent \qquad \qquad((4^41^4)^2(41)^2)^24^4((14)^2(1^44^4)^2)^21)^2(((4^41^4)^2(41)^2)^2(4^41)^2)^2(4^4\underline 1^4)^2\\
\indent \qquad\qquad ((41^4)^2((41)^2(4^41^4)^2)^2)^24\cdots \\
\indent \Delta^1(w)= (((41)^2(4^41^4)^2)^2(41^4)^2)^2(41)^2((4^41)^2((4^41^4)^2(41)^2)^2)^24^4((((14)^2(1^44^4)^2)^2(14^4)^2)^2(14)^2\\
\indent \qquad \qquad1^4(41)^24^4)^2((14)^2(1^44^4)^2)^2(1(4^41^4)^24(1^44^4)^2((14^4)^2((14)^2(1^44^4)^2)^2)^2)^21(((((4^41^4)^2\\
\indent \qquad \qquad(41)^2)^2(4^41)^2)^2(4^41^4)^24(1^44^4)^21)^2((4^41^4)^2(41)^2)^24^4((14)^2(1^44^4)^2)^21)^2 \cdots \\
\indent \Delta^2(w)= ((1^44^4)^2(14)^2)^21^4((41)^2(4^41^4)^2)^24(((1^44^4)^2(14)^2)^2(1^44)^2)^2(1^44^4)^2((14^4)^2((14)^2\\
\indent \qquad\qquad (1^44^4)^2)^2)^21((((4^41^4)^2(41)^2)^2(4^41)^2)^2(4^41^4)^24(1^44)^21)^2 \cdots \\
\indent \Delta^3(w)=(4^41^4)^24(1^44^4)^21((4^41^4)^2(41)^2)^24^4((14)^2(1^44^4)^2)^21(((4^41^4)^2(41)^2)^2(4^41)^2)^2 \cdots \\
\indent \Delta^4(w)=(4^41)^2(4^41^4)^24(1^44^4)^21((4^41^4)^2(41)^2)^2 \cdots \\
\indent \Delta^5(w)=(41)^2(4^41)^2(4^41^4)^2 \cdots \\
\indent \Delta^6(w)=1^4(41)^24^4 \cdots \\
\indent \Delta^7(w)=41^44 \cdots \\
\indent \Delta^8(w)=14 \cdots \\
\indent \Delta^9(w)=1 \cdots $ \\
$w$ has the factor $1^4(41^4)^2(41)^2$.

\item [(b)] If $p=1414441414$, then\\
\indent $\displaystyle \Delta^0(w)=(((1^44)^2((1^44^4)^2(14)^2)^2)^21^4(41)^24^4(14)^2)^2(1^44)^2((1^44^4)^21(4^41^4)^24\\
\indent \qquad \qquad (((1^44^4)^2(14)^2)^2(1^44)^2)^2)^2(1^44^4)^2 ((((14^4)^2((14)^2(1^44^4)^2)^2)^21(4^41^4)^24(1^44^4)^2)^2\\
\indent \qquad \qquad (14^4)^2(14)^2(1^44)^2(1^44^4)^2)^2((14^4)^2((14)^2(1^44^4)^2)^2)^21((4^41^4)^2(41)^2)^24^4((14)^2 \\
\indent \qquad \qquad  (1^44^4)^2)^21(4^41^4)^24(1^44^4)^2(14^4)^2(14)^2(1^44)^2(1^44^4)^2((14^4)^2((14)^2(1^44^4)^2)^2)^21\\
\indent \qquad \qquad (((4^41^4)^2(41)^2)^2(4^41)^2)^2 (4^4\underline 1^4)^2 (41^4)^2(41)^24^41\cdots \\
\indent  \Delta^1(w)=(((41)^2(4^41^4)^2)^2(41^4)^2)^2(41)^2((4^41)^2((4^41^4)^2(41)^2)^2)^24^4((((14)^2(1^44^4)^2)^2(14^4)^2)^2\\
\indent \qquad \qquad (14)^21^4(41)^24^4)^2((14)^2(1^44^4)^2)^21(4^41^4)^24(1^44^4)^2(14^4)^2(14)^21^4(41)^24^4\\
\indent \qquad \qquad((14)^2(1^44^4)^2)^21((4^41^4)^2(41)^2)^24^4 (14)^21^44\cdots \\
\indent \Delta^2(w)= ((1^44^4)^2(14)^2)^21^4((41)^2(4^41^4)^2)^24(((1^44^4)^2(14)^2)^2(1^44)^2)^2(1^44^4)^2(14^4)^2\\
\indent \qquad \qquad (14)^2(1^44)^2(1^44^4)^21(4^41^4)^24 1^44\cdots \\
\indent \Delta^3(w)=(4^41^4)^24(1^44^4)^21((4^41^4)^2(41)^2)^24^4(14)^21^4(41)^2(4^41)^2 4\cdots \\
\indent \Delta^4(w)=(4^41)^2(4^41^4)^2(41^4)^2(41)^2\cdots \\
\indent \Delta^5(w)=(41)^24^4(14)^21^4\cdots \\
\indent \Delta^6(w)=(1^44)^2\cdots \\
\indent \Delta^7(w)=(41)^2\cdots \\
\indent \Delta^8(w)=1^4\cdots \\
\indent \Delta^9(w)=4  \cdots $\\
$w$ has the factor $1^4(41^4)^2(41)^2$.

\item[(c)] If $p=141444144$, then\\
\indent $\displaystyle \Delta^0(w)=(((1^44)^2((1^44^4)^2(14)^2)^2)^21^4(41)^24^4(14)^2)^2(1^44)^2((1^44^4)^21(4^41^4)^24(((1^44^4)^2(14)^2)^2\\
\indent \qquad \qquad(1^44)^2)^2)^2(1^44^4)^2((((14^4)^2((14)^2(1^44^4)^2)^2)^21(4^41^4)^24(1^44^4))^2(14^4)^2(14)^2 \\
\indent \qquad \qquad (1^44)^2(1^44^4)^2)^2((14^4)^2((14)^2(1^44^4)^2)^2)^21(((4^41^4)^2(41)^2)^2(4^41)^2)^2(4^4\underline 1^4)^2\\
\indent \qquad \qquad ((((41^4)^2((41)^2(4^41^4)^2)^2)^2 4(1^44^4)^21(4^41^4))^2(41^4)^2(41)^2(4^41)^2(4^41^4)^2)^2\cdots \\
\indent  \Delta^1(w)=(((41)^2(4^41^4)^2)^2(41^4)^2)^2(41)^2((4^41)^2((4^41^4)^2(41)^2)^2)^24^4(((((14)^2(1^44^4)^2)^2\\
\indent \qquad \qquad (14^4)^2)^2(14)^21^4(41)^24^4)^2((14)^2(1^44^4)^2)^21((4^41^4)^2(41)^2)^24^4)^2 \cdots \\
\indent \Delta^2(w)= (((1^44^4)^2(14)^2)^21^4((41)^2(4^41^4)^2)^24((((1^44^4)^2(14)^2)^2(1^44)^2)^2(1^44^4)^21(4^41^4)^24)^2)^2\cdots \\
\indent \Delta^3(w)= ((4^41^4)^24(1^44^4)^21(((4^41^4)^2(41)^2)^2(4^41)^2)^2)^2 \cdots \\
\indent \Delta^4(w)= ((4^41)^2((4^41^4)^2(41)^2)^2)^2\cdots \\
\indent \Delta^5(w)= ((41)^2(4^41^4)^2)^2\cdots \\
\indent \Delta^6(w)= (1^44^4)^2\cdots \\
\indent  \Delta^7(w)=4^4\cdots \\
\indent  \Delta^8(w)=4\cdots $\\
$w$ has the factor $1^4(41^4)^2(41)^2$.
\end{itemize}
\end{itemize}

\item [Case (12)] If $p=1b1bbbb$, then\\
\indent $\displaystyle \Delta^0(w)= ((((1^bb)^\frac{b}{2}((1^bb^b)^\frac{b}{2}(1b)^\frac{b}{2})^\frac{b}{2})^\frac{b}{2}(1^b(b1)^\frac{b}{2}b^b(1b)^\frac{b}{2})^\frac{b}{4})^\frac{b}{2}((1^bb)^\frac{b}{2}(1^bb^b)^\frac{b}{2}(1b^b)^\frac{b}{2}(1b)^\frac{b}{2})^\frac{b}{4})^\frac{b}{2}\\
\indent \qquad \qquad ((\underline 1^bb)^\frac{b}{2}((1^bb^b)^\frac{b}{2}(1b)^\frac{b}{2})^\frac{b}{2})^\frac{b}{2}1^b((b1)^\frac{b}{2}(b^b1^b)^\frac{b}{2})^\frac{b}{2}b1^bb \cdots \\
\indent \Delta^1(w)= ((((b1)^\frac{b}{2}(b^b1^b)^\frac{b}{2})^\frac{b}{2}(b1^b)^\frac{b}{2})^\frac{b}{2}((b1)^\frac{b}{2}b^b(1b)^\frac{b}{2}1^b)^\frac{b}{4})^\frac{b}{2}((b1)^\frac{b}{2}(b^b1^b)^\frac{b}{2})^\frac{b}{2}b(1^bb^b)^\frac{b}{2}1b\cdots \\
\indent \Delta^2(w)= (((1^bb^b)^\frac{b}{2}(1b)^\frac{b}{2})^\frac{b}{2}(1^bb)^\frac{b}{2})^\frac{b}{2}(1^bb^b)^\frac{b}{2}1b^b1\cdots \\
\indent \Delta^3(w)= ((b^b1^b)^\frac{b}{2}(b1)^\frac{b}{2})^\frac{b}{2}b^b1b\cdots \\
\indent \Delta^4(w)= (b^b1^b)^\frac{b}{2}b1\cdots \\
\indent \Delta^5(w)= b^b1\cdots \\
\indent \Delta^6(w)= b\cdots $\\
$w$ has the factor $((1^bb)^\frac{b}{2}((1^bb^b)^\frac{b}{2}(1b)^\frac{b}{2})^\frac{b}{2})^\frac{b}{2}1^b(b1)^\frac{b}{2}(b^b1^b)^\frac{b}{2}$. 

\item [Case (13)] If $p=1bb111$, then from Case (6), we have\\
\indent $\displaystyle \Delta^0(w)= (1^bb^b)((1b^b)^\frac{b}{2}((1b)^\frac{b}{2}(\underline 1^bb^b)^\frac{b}{2})^\frac{b}{2})^\frac{b}{2}1(b^b1^b)^\frac{b}{2}b1\cdots \\$
$w$ has the factor $(1^bb^b)^\frac{b}{2}1(b^b1^b)^\frac{b}{2}$.

\item [Case (14)] If $p=1bb11b$, using Case (7), we get \\
\indent $\displaystyle \Delta^0(w)=((1^bb^b)^\frac{b}{2}(1b^b)^\frac{b}{2}(1b)^\frac{b}{2}(\underline 1^bb)^\frac{b}{2})^\frac{b}{4} ((1^bb^b)^\frac{b}{2}1(b^b1^b)^\frac{b}{2}b)^\frac{b}{4}((1^bb^b)^\frac{b}{2}(1b)^\frac{b}{2})^\frac{b}{2}1^bb1 \cdots \\$
$w$ has the factor $(1^bb)^\frac{b}{2}$. 

\item [Case (15)] If $p=1bb1b$, using Case (8), we get\\
\indent $\displaystyle \Delta^0(w)=((1^bb^b)^\frac{b}{2}1(b^b\underline 1^b)^\frac{b}{2}b)^\frac{b}{4}((1^bb^b)^\frac{b}{2}(1b)^\frac{b}{2})^\frac{b}{2}1^bb \cdots \\$
$w$ has the factor $1^bb1^b$.

\item [Case (16)] If $p=1bbb$, then\\
\indent $\displaystyle \Delta^0(w)=((1^bb^b)^\frac{b}{2}(1b)^\frac{b}{2})^\frac{b}{2}\underline 1^bb1 \cdots \\
\indent \Delta^1(w)=(b^b1^b)^\frac{b}{2}b1 \cdots \\
\indent \Delta^2(w)= b^b1\cdots \\
\indent \Delta^3(w)=b \cdots $\\
$w$ has the factor $1^bb1$.
\end{itemize}

\section{Details of the proof of Theorem \ref{thm122n}}

\begin{itemize}
\item [Case (17)] If $p=1b11b1$, then\\
\indent $\displaystyle  \Delta^0(w)= 1^b(b1)^{2n+1}(b^b(1b)^{2n+1}1^b(b1)^{2n+1})^n(b^b1)^{2n+1}(b^b\underline 1^b)^{2n+1}b1^bb\cdots \\
\indent \Delta^1(w)= (b1^b)^{2n+1}(b1)^{2n+1}b^b1b\cdots \\
\indent \Delta^2(w)=(1b)^{2n+1}1^bb1\cdots \\
\indent \Delta^3(w)=1^bb1\cdots \\
\indent \Delta^4(w)=b1\cdots \\
\indent \Delta^5(w)=1\cdots $\\
$w$ has the factor $1^bb1^b$.

\item [Case (18)] If $p=1b11bb$, then\\
\indent $\displaystyle \Delta^0(w)=(1^b(b1)^{2n+1}(b^b(1b)^{2n+1}1^b(b1)^{2n+1})^n((b^b1)^{2n+1}((b^b\underline 1^b)^{2n+1}(b1)^{2n+1})^{2n+1})^{2n+1})^{2n+1}\\
\indent \qquad \qquad b^b(1b)^{2n+1}(1^bb^b)^{2n+1}\cdots \\
\indent \Delta^1(w)=((b1^b)^{2n+1}((b1)^{2n+1}(b^b1^b)^{2n+1})^{2n+1})^{2n+1}b(1^bb^b)^{2n+1}1b\cdots \\
\indent \Delta^2(w)=((1b)^{2n+1}(1^bb^b)^{2n+1})^{2n+1}1b^b1\cdots \\
\indent \Delta^3(w)=(1^bb^b)^{2n+1}1b\cdots \\
\indent \Delta^4(w)= b^b1\cdots \\
\indent \Delta^5(w)= b\cdots $\\
$w$ has the factor $1^bb1b^b(1b)1^b$ if $n=0$ and $1^b(b1)^{2n+1}(b^b1^b)^{2n+1}$ in $((b^b1^b)^{2n+1}(b1)^{2n+1})^{2n+1}$ if $n\neq 0$.

\item [Case (19)] If $p=1b1b11b11$, there are $2$ subcases to consider.
\begin{itemize}
\item [(i)] If $b=2$, then\\
\indent $\displaystyle \Delta^0(w)=1121122122112122121122122121121122122112122122\underline112 112122112112212\cdots \\
\indent \Delta^1(w)=21221221121122121121221221121221211221221\cdots \\
\indent \Delta^2(w)=112122122112112122112112212\cdots \\
\indent \Delta^3(w)=21121221211221221\cdots \\
\indent \Delta^4(w)=12112112212\cdots \\
\indent \Delta^5(w)=1121221\cdots \\
\indent \Delta^6(w)=2112\cdots \\
\indent \Delta^7(w)=12\cdots \\
\indent \Delta^8(w)=1\cdots $\\
$w$ has the factor $1121121$. 

\item [(ii)] If $b\neq 2$, then $n\neq 0$ and \\
\indent $\displaystyle \Delta^0(w)= (1^bb)^{2n+1}(1^bb^b)^{2n+1}1((b^b1^b)^{2n+1}b(1^bb^b)^{2n+1}1)^n((b^b1^b)^{2n+1}(b1)^{2n+1})^{2n+1}b^b(1b)^{2n+1}\\
\indent \qquad \qquad(1^b(b1)^{2n+1}b^b(1b)^{2n+1})^n(1^bb)^{2n+1}(1^bb^b)^{2n+1}1 ((b^b1^b)^{2n+1}b(1^bb^b)^{2n+1}1)^n)^n\\
\indent \qquad \qquad ((b^b1^b)^{2n+1}(b1)^{2n+1})^{2n+1}b^b(((1b)^{2n+1}(1^bb^b)^{2n+1})^{2n+1}1((b^b1^b)^{2n+1}(b1)^{2n+1})^{2n+1}b^b)^n\\
\indent \qquad \qquad (((1b)^{2n+1}(1^bb^b)^{2n+1})^{2n+1}(1b^b)^{2n+1})^{2n+1}((1b)^{2n+1}(1^bb)^{2n+1} ((1^bb^b)^{2n+1}(1b^b)^{2n+1}(1b)^{2n+1}\\
\indent \qquad \qquad (1^bb)^{2n+1})^n(1^bb^b)^{2n+1}1((b^b1^b)^{2n+1}b(1^bb^b)^{2n+1}1)^n((b^b1^b)^{2n+1}(b1)^{2n+1})^{2n+1}b^b\\
\indent \qquad \qquad(((1b)^{2n+1}(1^bb^b)^{2n+1})^{2n+1}1((b^b1^b)^{2n+1}(b1)^{2n+1})^{2n+1}b^b)^n(((1b)^{2n+1}(1^bb^b)^{2n+1})^{2n+1}\\
\indent \qquad \qquad(1b^b)^{2n+1})^{2n+1})^n(1b)^{2n+1}(1^bb)^{2n+1}((1^bb^b)^{2n+1}(1b^b)^{2n+1}(1b)^{2n+1}(1^bb)^{2n+1})^n(1^bb^b)^{2n+1}1\\
\indent \qquad \qquad ((b^b1^b)^{2n+1}b(1^bb^b)^{2n+1}1)^n(((b^b1^b)^{2n+1}(b1)^{2n+1})^{2n+1}(b^b1)^{2n+1})^{2n+1}((b^b1^b)^{2n+1}b\\
\indent \qquad \qquad((1^bb^b)^{2n+1}1(b^b1^b)^{2n+1}b)^n((1^bb^b)^{2n+1}(1b)^{2n+1})^{2n+1}( 1^bb)^{2n+1})^{2n+1}(1^bb^b)^{2n+1}1\\
\indent \qquad \qquad ((b^b1^b)^{2n+1}b(1^bb^b)^{2n+1}1)^n((b^b1^b)^{2n+1}(b1)^{2n+1})^{2n+1}(b^b1)^{2n+1})^{2n+1})^n ((b^b\underline1^b)^{2n+1}\\
\indent \qquad \qquad (b1^b)^{2n+1}((b1)^{2n+1}(b^b1)^{2n+1}(b^b1^b)^{2n+1}(b1^b)^{2n+1})^n((b1)^{2n+1}b^b((1b)^{2n+1}1^b(b1)^{2n+1}b^b)^n\cdots\\
\indent \Delta^1(w)=(b1)^{2n+1}(b^b1)^{2n+1}((b^b1^b)^{2n+1}(b1^b)^{2n+1}(b1)^{2n+1}(b^b1)^{2n+1})^n(b^b1^b)^{2n+1}b((1^bb^b)^{2n+1}1\\
\indent \qquad \qquad (b^b1^b)^{2n+1}b)^n ((1^bb^b)^{2n+1}(1b)^{2n+1})^{2n+1}(1^b(b1)^{2n+1}(b^b(1b)^{2n+1}1^b(b1)^{2n+1})^n(b^b1)^{2n+1}\\
\indent \qquad \qquad (b^b1^b)^{2n+1}b((1^bb^b)^{2n+1}1(b^b1^b)^{2n+1}b)^n((1^bb^b)^{2n+1}(1b)^{2n+1})^{2n+1})^n  1^b(b1)^{2n+1}(b^b(1b)^{2n+1}1^b\\
\indent \qquad \qquad (b1)^{2n+1})^n((b^b1)^{2n+1}((b^b1^b)^{2n+1}(b1)^{2n+1})^{2n+1})^{2n+1} (b^b(1b)^{2n+1}(1^b(b1)^{2n+1}b^b(1b)^{2n+1})^n\\
\indent \qquad \qquad ((1^bb)^{2n+1}((1^bb^b)^{2n+1}(1b)^{2n+1})^{2n+1})^{2n+1} 1^b(b1)^{2n+1}(b^b(1b)^{2n+1}1^b(b1)^{2n+1})^n((b^b1)^{2n+1}\\
\indent \qquad \qquad ((b^b1^b)^{2n+1}(b1)^{2n+1})^{2n+1})^{2n+1})^nb^b((1b)^{2n+1}(1^bb^b)^{2n+1})^{2n+1}1b^b1\cdots \\
\indent \Delta^2(w)=1^b(b1)^{2n+1}(b^b(1b)^{2n+1}1^b(b1)^{2n+1})^n(b^b1)^{2n+1}(b^b1^b)^{2n+1}((b1^b)^{2n+1}(b1)^{2n+1}(b^b1)^{2n+1}\\
\indent \qquad \qquad (b^b1^b)^{2n+1})^n((b1^b)^{2n+1}((b1)^{2n+1}(b^b1^b)^{2n+1})^{2n+1})^{2n+1}b(1^bb^b)^{2n+1}1b\cdots \\
\indent \Delta^3(w)=(b1^b)^{2n+1}(b1)^{2n+1}b^b((1b)^{2n+1}1^b(b1)^{2n+1}b^b)^n((1b)^{2n+1}(1^bb^b)^{2n+1})^{2n+1}1b^b1\cdots \\
\indent \Delta^4(w)=(1b)^{2n+1}(1^bb)^{2n+1}(1^bb^b)^{2n+1}1b\cdots \\
\indent \Delta^5(w)= 1^b(b1)^{2n+1}b1\cdots \\
\indent \Delta^6(w)=b1^bb\cdots \\
\indent \Delta^7(w)=1b\cdots \\
\indent \Delta^8(w)= 1\cdots $\\
$w$ has the factor $1^b(b1^b)^{2n+1}(b1)^{2n+1}$.

\end{itemize}

\item [Case (20)] If $p=1b1b11b1b$, there are $2$ subcases to consider.
\begin{itemize}
\item [(i)] If $b=2$, then\\
\indent $\displaystyle \Delta^0(w)=112112212211212212112212212\underline 1121122121121221221121121221121122121121221221121 \cdots \\
\indent \Delta^1(w)=21221221121122121121221121121221211221221121121221\cdots \\
\indent \Delta^2(w)=112122122112112212112112212212112\cdots \\
\indent \Delta^3(w)=211212212211212212112\cdots \\
\indent \Delta^4(w)=12112122112112\cdots \\
\indent \Delta^5(w)=112112212\cdots \\
\indent \Delta^6(w)=21221\cdots \\
\indent \Delta^7(w)=112\cdots \\
\indent \Delta^8(w)=2\cdots $\\
$w$ has the factor $1121122121$.
\item [(ii)] If $b \neq 2$, then $n\neq 0$ and \\
\indent $\displaystyle \Delta^0(w)=(1^bb)^{2n+1}(1^bb^b)^{2n+1}1((b^b1^b)^{2n+1}b(1^bb^b)^{2n+1}1)^n(((b^b1^b)^{2n+1}(b1)^{2n+1})^{2n+1}b^b(1b)^{2n+1}(1^b\\
\indent \qquad \qquad (b1)^{2n+1}b^b(1b)^{2n+1})^n(1^bb)^{2n+1}(1^bb^b)^{2n+1}1((b^b1^b)^{2n+1}b (1^bb^b)^{2n+1}1)^n)^n\\
\indent \qquad \qquad ((b^b1^b)^{2n+1}(b1)^{2n+1})^{2n+1}b^b(((1b)^{2n+1}(1^bb^b)^{2n+1})^{2n+1}1((b^b1^b)^{2n+1}(b1)^{2n+1})^{2n+1}b^b)^n\\
\indent \qquad \qquad ((((1b)^{2n+1}(1^bb^b)^{2n+1})^{2n+1}(1b^b)^{2n+1})^{2n+1}(1b)^{2n+1} (((\underline 1^bb)^{2n+1}((1^bb^b)^{2n+1}\\
\indent \qquad \qquad (1b)^{2n+1})^{2n+1})^{2n+1}1^b \cdots \\
\indent \Delta^1(w)=(b1)^{2n+1}(b^b1)^{2n+1}((b^b1^b)^{2n+1}(b1^b)^{2n+1}(b1)^{2n+1}(b^b1)^{2n+1})^n(b^b1^b)^{2n+1}b((1^bb^b)^{2n+1}1\\
\indent \qquad \qquad (b^b1^b)^{2n+1}b)^n(((1^bb^b)^{2n+1}(1b)^{2n+1})^{2n+1}1^b(((b1)^{2n+1}(b^b1^b)^{2n+1})^{2n+1}b((1^bb^b)^{2n+1} \\
\indent \qquad \qquad (1b)^{2n+1})^{2n+1}1^b)^n\cdots \\
\indent \Delta^2(w)= 1^b(b1)^{2n+1}(b^b(1b)^{2n+1}1^b(b1)^{2n+1})^n(b^b1)^{2n+1}((b^b1^b)^{2n+1}b((1^bb^b)^{2n+1}1(b^b1^b)^{2n+1}b)^n\\
\indent \qquad \qquad ((1^bb^b)^{2n+1}(1b)^{2n+1})^{2n+1}1^b(b1)^{2n+1}(b^b(1b)^{2n+1}1^b(b1)^{2n+1})^n(b^b1)^{2n+1})^n\cdots \\
\indent \Delta^3(w)=(b1^b)^{2n+1}(b1)^{2n+1}((b^b1)^{2n+1}(b^b1^b)^{2n+1}(b1^b)^{2n+1}(b1)^{2n+1})^n((b^b1)^{2n+1}\\
\indent \qquad \qquad ((b^b1^b)^{2n+1}(b1)^{2n+1})^{2n+1})^{2n+1}b^b(1b)^{2n+1}1^bb\cdots \\
\indent \Delta^4(w)= (1b)^{2n+1}1^b((b1)^{2n+1}b^b(1b)^{2n+1}1^b)^n((b1)^{2n+1}(b^b1^b)^{2n+1})^{2n+1}b1^bb\cdots \\
\indent \Delta^5(w)=(1^bb)^{2n+1}(1^bb^b)^{2n+1}1b\cdots \\
\indent \Delta^6(w)=(b1)^{2n+1}b^b1\cdots \\
\indent \Delta^7(w)=1^bb\cdots \\
\indent \Delta^8(w)=b\cdots $\\
$w$ has the factor $(1^bb)^{2n+1}(1^bb^b)^{2n+1}(1b)^{2n+1}$.
\end{itemize}

\item [Case (21)] If $p=1b1b11bb$, there are $2$ subcases to consider.
\begin{itemize}
\item [(i)] If $b=2$, there are again 3 subcases:
\begin{itemize}
\item [(a)] If $p=1212112211$, then\\
\indent $\displaystyle \Delta^0(w)=11211221221121221211221221211212211211221221211221221121221211 \\
\indent \qquad \qquad 21122122121122122112122122\underline 1121122121121221\cdots\\
\indent \Delta^1(w)=21221221121122121121122122121122122112112122121122122112122122122\\
\indent \qquad \qquad 112112\cdots\\
\indent \Delta^2(w)=11212212211212212112212212112112212211212212\cdots\\
\indent \Delta^3(w)= 21121221121122121121221221121\cdots\\
\indent \Delta^4(w)=1211221221121121221 \cdots\\
\indent \Delta^5(w)=112212212112\cdots\\
\indent \Delta^6(w)=2212112\cdots\\
\indent  \Delta^7(w)=2112\cdots\\
\indent \Delta^8(w)=12\cdots\\
\indent \Delta^9(w)=1\cdots$\\
$w$ has the factor $11211221211$. 

\item [(b)] If $p=1212112212$, then\\
\indent $\displaystyle\Delta^0(w)= 112112212211212212112212212112122112112212212112212211212212112\\
\indent \qquad \qquad 1122122112122122 \underline 11211212211211221221211211221211\cdots\\
\indent \Delta^1(w)=  212212211211221211211221221211221221121121221\\ 
\indent \qquad \qquad 2211212212112212212112122112112212212112\cdots\\
\indent \Delta^2(w)=11212212211212212112212212112122112112212112112212212112
\cdots\\
\indent  \Delta^3(w)=2112122112112212112112212211212212112\cdots\\
\indent \Delta^4(w)=121122122112122122112112\cdots\\
\indent \Delta^5(w)=112212211212212\cdots\\
\indent \Delta^6(w)=221221121\cdots\\
\indent \Delta^7(w)=21221\cdots\\
\indent \Delta^8(w)=112\cdots\\
\indent \Delta^9(w)=2 \cdots$\\
$w$ has the factor $1121121$.

\item [(c)] If $p=121211222$, then\\
\indent $\displaystyle\Delta^0(w)= 1121122122112122121122122121121221121122122121122122112122122\underline 11211\\
\indent \qquad\qquad 21221211221221121\cdots\\
\indent \Delta^1(w)=21221221121122 1211211221221211221221121221211211221221\cdots\\
\indent \Delta^2(w)=11212212211212212112212211211212212\cdots\\
\indent \Delta^3(w)=21121221121122122121121\cdots\\
\indent \Delta^4(w)=121122122121121\cdots\\
\indent \Delta^5(w)=1122121121\cdots\\
\indent \Delta^6(w)=221121\cdots\\
\indent \Delta^7(w)=221 \cdots\\
\indent \Delta^8(w)=2 \cdots$\\
$w$ has the factor $1121121$.
\end{itemize}

\item [(ii)] If $b\neq 2$, then $n\neq 0$ and \\
\indent $\displaystyle \Delta^0(w)=((1^bb)^{2n+1}(1^bb^b)^{2n+1}1((b^b1^b)^{2n+1}b(1^bb^b)^{2n+1}1)^n(((b^b1^b)^{2n+1}(b1)^{2n+1})^{2n+1}\\
\indent \qquad \qquad b^b(1b)^{2n+1}(1^b(b1)^{2n+1}b^b(1b)^{2n+1})^n(1^bb)^{2n+1}(1^bb^b)^{2n+1}1((b^b1^b)^{2n+1}b(1^bb^b)^{2n+1}1)^n)^n\\
\indent \qquad \qquad ((b^b1^b)^{2n+1}(b1)^{2n+1})^{2n+1} b^b(((1b)^{2n+1}(1^bb^b)^{2n+1})^{2n+1}1((b^b1^b)^{2n+1}(b1)^{2n+1})^{2n+1}\\
\indent \qquad \qquad b^b)^n(((1b)^{2n+1}(1^bb^b)^{2n+1})^{2n+1}(1b^b)^{2n+1})^{2n+1}(1b)^{2n+1}1^b((b1)^{2n+1}b^b(1b)^{2n+1}1^b)^n(((b1)^{2n+1}\\
\indent \qquad \qquad (b^b\underline1^b)^{2n+1})^{2n+1}(b1^b)^{2n+1})^{2n+1}(b1)^{2n+1}b^b((1b)^{2n+1}1^b(b1)^{2n+1}b^b)^n \cdots \\
\indent \Delta^1(w)=((b1)^{2n+1}(b^b1)^{2n+1}((b^b1^b)^{2n+1}(b1^b)^{2n+1}(b1)^{2n+1}(b^b1)^{2n+1})^n(b^b1^b)^{2n+1}b((1^bb^b)^{2n+1}\\
\indent \qquad \qquad 1(b^b1^b)^{2n+1}b)^n (((1^bb^b)^{2n+1}(1b)^{2n+1})^{2n+1}(1^bb)^{2n+1})^{2n+1})^{2n+1}\cdots\\
\indent \Delta^2(w)=(1^b(b1)^{2n+1}(b^b(1b)^{2n+1}1^b(b1)^{2n+1})^n((b^b1)^{2n+1}((b^b1^b)^{2n+1}(b1)^{2n+1})^{2n+1})^{2n+1})^{2n+1}\cdots \\
\indent \Delta^3(w)=((b1^b)^{2n+1}((b1)^{2n+1}(b^b1^b)^{2n+1})^{2n+1})^{2n+1}\cdots\\
\indent \Delta^4(w)=((1b)^{2n+1}(1^bb^b)^{2n+1})^{2n+1}\cdots\\
\indent \Delta^5(w)=(1^bb^b)^{2n+1}\cdots\\
\indent \Delta^6(w)= b^b\cdots\\
\indent \Delta^7(w)=b\cdots$\\
$w$ has the factor $1^b(b1^b)^{2n+1}(b1)^{2n+1}$.
\end{itemize}

\item [Case (22)] If $p=1b1b1b11$, there are $2$ subcases to consider.
\begin{itemize}
\item [(i)] If $b=2$, then\\
\indent $\displaystyle \Delta^0(w)= 1121122122121122122112122122\underline 1121122121121221\cdots\\
\indent \Delta^1(w)=2122121122122112122122112112\cdots\\
\indent \Delta^2(w)=112112212211212212\cdots\\
\indent \Delta^3(w)=21221221121\cdots\\
\indent \Delta^4(w)=1121221\cdots\\
\indent \Delta^5(w)=2112\cdots\\
\indent \Delta^6(w)=12\cdots\\
\indent \Delta^7(w)=1\cdots$\\
$w$ has the factor $1121122121$.

\item[(ii)] If $b\neq 2$, then $n\neq 0$ and \\
\indent $\displaystyle\Delta^0(w)= (1^bb)^{2n+1}(1^bb^b)^{2n+1}((1b^b)^{2n+1}(1b)^{2n+1}(1^bb)^{2n+1}(1^bb^b)^{2n+1})^n1(b^b1^b)^{2n+1}(b(1^bb^b)^{2n+1}\\
\indent \qquad \qquad 1(b^b\underline 1^b)^{2n+1})^n(b1^b)^{2n+1}((b1)^{2n+1}b^b((1b)^{2n+1}1^b(b1)^{2n+1}b^b)^n((1b)^{2n+1}(1^bb^b)^{2n+1})^{2n+1}\\
\indent \qquad \qquad 1(b^b1^b)^{2n+1}(b(1^bb^b)^{2n+1}1(b^b1^b)^{2n+1})^n(b1^b)^{2n+1})^n\cdots \\
\indent \Delta^1(w)=(b1)^{2n+1}b^b((1b)^{2n+1}1^b(b1)^{2n+1}b^b)^n(1b^b)^{2n+1}(1b)^{2n+1}((1^bb)^{2n+1}(1^bb^b)^{2n+1}(1b^b)^{2n+1}\\
\indent \qquad \qquad (1b)^{2n+1})^n((1^bb)^{2n+1}((1^bb^b)^{2n+1}(1b)^{2n+1})^{2n+1})^{2n+1}1^b(b1)^{2n+1}b^b1\cdots\\
\indent \Delta^2(w)=(1^bb)^{2n+1}(1b)^{2n+1}1^b((b1)^{2n+1}b^b(1b)^{2n+1}1^b)^n((b1)^{2n+1}(b^b1^b)^{2n+1})^{2n+1}b1^bb\cdots\\
\indent \Delta^3(w)=(b1)^{2n+1}(1^bb)^{2n+1}(1^bb^b)^{2n+1}1b\cdots\\
\indent \Delta^4(w)=1^b(b1)^{2n+1}b^b1\cdots\\
\indent \Delta^5(w)=b1^bb\cdots\\
\indent \Delta^6(w)=1b\cdots\\
\indent \Delta^7(w)=1\cdots$\\
$w$ has the factor $1^b(b1^b)^{2n+1}(b1)^{2n+1}$.
\end{itemize}

\item [Case (23)] If $p=1b1b1b1b$, there are $2$ cases to consider.
\begin{itemize}
\item [(i)] If $b=2$, then\\
\indent $\displaystyle \Delta^0(w)=112112212212112212211212212112112212112122\underline 11211212212\cdots\\
\indent \Delta^1(w)=21221211221221121121221121122121121\cdots\\
\indent \Delta^2(w)=11211221221211221221121\cdots\\
\indent \Delta^3(w)=21221211221221\cdots\\
\indent \Delta^4(w)=112112212\cdots\\
\indent \Delta^5(w)=21221\cdots\\
\indent \Delta^6(w)=112\cdots\\
\indent \Delta^7(w)=2\cdots$\\
$w$ has the factor $1121121$.

\item [(ii)] If $b\neq 2$, then\\
\indent $\displaystyle \Delta^0(w)=(1^bb)^{2n+1}(1^bb^b)^{2n+1}((1b^b)^{2n+1}(1b)^{2n+1}(1^bb)^{2n+1}(1^bb^b)^{2n+1})^n((1b^b)^{2n+1}((1b)^{2n+1}\\
\indent \qquad \qquad (1^bb^b)^{2n+1})^{2n+1})^{2n+1}(1((b^b1^b)^{2n+1}(b1)^{2n+1})^{2n+1}(b^b((1b)^{2n+1}(1^bb^b)^{2n+1})^{2n+1}1((b^b1^b)^{2n+1}\\
\indent \qquad \qquad (b1)^{2n+1})^{2n+1})^nb^b(1b)^{2n+1}(1^b(b1)^{2n+1}b^b(1b)^{2n+1})^n( 1^bb)^{2n+1}(1^bb^b)^{2n+1}((1b^b)^{2n+1}\\
\indent \qquad \qquad(1b)^{2n+1}(1^bb)^{2n+1}(1^bb^b)^{2n+1})^n((1b^b)^{2n+1}((1b)^{2n+1}(1^bb^b)^{2n+1})^{2n+1})^{2n+1})^n 1((b^b1^b)^{2n+1}\\
\indent \qquad \qquad(b1)^{2n+1})^{2n+1}(b^b((1b)^{2n+1}(1^bb^b)^{2n+1})^{2n+1}1((b^b1^b)^{2n+1}(b1)^{2n+1})^{2n+1})^nb^b(1b)^{2n+1}\\
\indent \qquad \qquad (1^b(b1)^{2n+1}b^b(1b)^{2n+1})^n((\underline1^bb)^{2n+1}((1^bb^b)^{2n+1}(1b)^{2n+1})^{2n+1})^{2n+1} \cdots\\
\indent \Delta^1(w)=(b1)^{2n+1}b^b((1b)^{2n+1}1^b(b1)^{2n+1}b^b)^n((1b)^{2n+1}(1^bb^b)^{2n+1})^{2n+1}(1(b^b1^b)^{2n+1}(b(1^bb^b)^{2n+1}1\\
\indent \qquad \qquad (b^b1^b)^{2n+1})^n(b1^b)^{2n+1}(b1)^{2n+1}b^b((1b)^{2n+1}1^b(b1)^{2n+1}b^b)^n((1b)^{2n+1}(1^bb^b)^{2n+1})^{2n+1})^n\\
\indent \qquad \qquad 1(b^b1^b)^{2n+1}(b(1^bb^b)^{2n+1}1(b^b1^b)^{2n+1})^n((b1^b)^{2n+1}((b1)^{2n+1}(b^b1^b)^{2n+1})^{2n+1})^{2n+1} \cdots \\
\indent \Delta^2(w)=(1^bb)^{2n+1}(1^bb^b)^{2n+1}((1b^b)^{2n+1}(1b)^{2n+1}(1^bb)^{2n+1}(1^bb^b)^{2n+1})^n((1b^b)^{2n+1}((1b)^{2n+1}\\
\indent \qquad \qquad (1^bb^b)^{2n+1})^{2n+1})^{2n+1}1(b^b1^b)^{2n+1}b1\cdots\\
\indent \Delta^3(w)=(b1)^{2n+1}b^b((1b)^{2n+1}1^b(b1)^{2n+1}b^b)^n((1b)^{2n+1}(1^bb^b)^{2n+1})^{2n+1}1b^b1\cdots\\
\indent \Delta^4(w)=(1^bb)^{2n+1}(1^bb^b)^{2n+1}1b\cdots\\
\indent \Delta^5(w)= (b1)^{2n+1}b^b1\cdots\\
\indent \Delta^6(w)=1^bb\cdots\\
\indent \Delta^7(w)=b\cdots$\\
$w$ has the factor $(1^bb)^{2n+1}(1^bb^b)^{2n+1}(1b)^{2n+1}$.
\end{itemize}

\item [Case (24)] If $p=1b1b1bb$, there are $2$ subcases to consider. 
\begin{itemize}
\item [(i)] If $b=2$, then\\
\indent $\displaystyle \Delta^0(w)=1121122122121122122\underline 1121121221211221221\cdots\\
\indent \Delta^1(w)=212212112212212112112212\cdots\\
\indent \Delta^2(w)=1121122121121221\cdots\\
\indent \Delta^3(w)=2122112112\cdots\\
\indent \Delta^4(w)=112212\cdots\\
\indent \Delta^5(w)=221\cdots\\
\indent \Delta^6(w)=2\cdots$\\
$w$ has the factor $1121121$.

\item [(ii)] If $b\neq 2$, then\\
\indent $\displaystyle \Delta^0(w)=((1^bb)^{2n+1}(1^bb^b)^{2n+1}((1b^b)^{2n+1}(1b)^{2n+1}(1^bb)^{2n+1}(1^bb^b)^{2n+1})^n((1b^b)^{2n+1}((1b)^{2n+1}\\
\indent \qquad \qquad (1^bb^b)^{2n+1})^{2n+1})^{2n+1}1(b^b1^b)^{2n+1}(b(1^bb^b)^{2n+1}1(b^b\underline 1^b)^{2n+1})^n ((b1^b)^{2n+1}((b1)^{2n+1}\cdots \\
\indent \Delta^1(w)=((b1)^{2n+1}b^b((1b)^{2n+1}1^b(b1)^{2n+1}b^b)^n(((1b)^{2n+1}(1^bb^b)^{2n+1})^{2n+1}(1b^b)^{2n+1})^{2n+1})^{2n+1}\cdots\\
\indent\Delta^2(w)= ((1^bb)^{2n+1}((1^bb^b)^{2n+1}(1b)^{2n+1})^{2n+1})^{2n+1}\cdots\\
\indent \Delta^3(w)=((b1)^{2n+1}(b^b1^b)^{2n+1})^{2n+1}\cdots\\
\indent \Delta^4(w)=(1^bb^b)^{2n+1}\cdots\\
\indent \Delta^5(w)=b^b\cdots\\
\indent \Delta^6(w)=b\cdots$\\
$w$ has the factor $1^b(b1^b)^{2n+1}(b)^{2n+1}$.
\end{itemize}

\item [Case (25)] If $p=1b1bb111$, there are $2$ subcases to consider.
\begin{itemize}
\item [(i)] If  $b=2$, there are again 2 subcases to consider.
\begin{itemize}
\item [(a)] If $p=121221111$, then\\
\indent $\displaystyle \Delta^0(w)=1121122121121221121122122121122122112122122\underline 112112122\cdots\\
\indent \Delta^1(w)=212211211221221211221221121221211211221221\cdots\\
\indent \Delta^2(w)=112212212112212211211212212\cdots\\
\indent \Delta^3(w)=22121122122121121\cdots\\
\indent \Delta^4(w)=21122121121\cdots\\
\indent \Delta^5(w)=1221121\cdots\\
\indent \Delta^6(w)=1221\cdots\\
\indent \Delta^7(w)=12\cdots\\
\indent \Delta^8(w)=1\cdots$\\
$w$ has the factor $1121121$.

\item [(b)] If $p=121221112$, then\\
\indent $\displaystyle \Delta^0(w)=1121122121121221121122122121122122112122122112112212112122122\underline 1121121\cdots\\
\indent \Delta^1(w)=212211211221221211221221121221221121121221211221221\cdots\\
\indent \Delta^2(w)=112212212112212211212212112112212\cdots\\
\indent \Delta^3(w)=221211221221121121221\cdots\\
\indent \Delta^4(w)=2112212212112\cdots\\
\indent \Delta^5(w)=12212112\cdots\\
\indent \Delta^6(w)=12112\cdots\\
\indent \Delta^7(w)=112\cdots\\
\indent \Delta^8(w)=2\cdots$\\
$w$ has the factor $1121121$.
\end{itemize}

\item [(ii)] If $b \neq 2$, then\\
\indent $\displaystyle \Delta^0(w)=((1^bb)^{2n+1}((1^bb^b)^{2n+1}(1b)^{2n+1})^{2n+1})^{2n+1}(1^b((b1)^{2n+1}(b^b1^b)^{2n+1})^{2n+1}(b((1^bb^b)^{2n+1}\\
\indent \qquad \qquad (1b)^{2n+1})^{2n+1}1^b((b1)^{2n+1}(b^b1^b)^{2n+1})^{2n+1})^nb(1^bb^b)^{2n+1}(1(b^b1^b)^{2n+1}b(1^bb^b)^{2n+1})^n\\
\indent \qquad \qquad((1b^b)^{2n+1}((1b)^{2n+1}(1^bb^b)^{2n+1})^{2n+1})^{2n+1}(1(b^b1^b)^{2n+1}(b(1^bb^b)^{2n+1}1(b^b\underline 1^b)^{2n+1})^n\\
\indent \qquad \qquad((b1^b)^{2n+1}((b1)^{2n+1}(b^b1^b)^{2n+1})^{2n+1})^{2n+1}\cdots \\
\indent \Delta^1(w)= ((b1)^{2n+1}(b^b1^b)^{2n+1})^{2n+1}(b(1^bb^b)^{2n+1}(1(b^b1^b)^{2n+1}b(1^bb^b)^{2n+1})^n((1b^b)^{2n+1}((1b)^{2n+1}\\
\indent \qquad \qquad(1^bb^b)^{2n+1})^{2n+1})^{2n+1})^{2n+1}1((b^b1^b)^{2n+1}(b1)^{2n+1})^{2n+1}b^b1b\cdots\\
\indent \Delta^2(w)=(1^bb^b)^{2n+1}((1b^b)^{2n+1}((1b)^{2n+1}(1^bb^b)^{2n+1})^{2n+1})^{2n+1}1(b^b1^b)^{2n+1}b1\cdots\\
\indent \Delta^3(w)=b^b((1b)^{2n+1}(1^bb^b)^{2n+1})^{2n+1}1b^b1\cdots\\
\indent \Delta^4(w)=b(1^bb^b)^{2n+1}1b\cdots\\
\indent \Delta^5(w)=1b^b1\cdots\\
\indent \Delta^6(w)=1b\cdots\\
\indent \Delta^7(w)=1\cdots$\\
$w$ has the factor $1^b(b1^b)^{2n+1}(b1)^{2n+1}$.

\end{itemize}

\item [Case (26)] If $p=1b1bb11b$, there are $2$ subcases to consider.
\begin{itemize}
\item [(i)] If $b=2$, then\\
\indent $\displaystyle \Delta^0(w)= 112112212112122112112212212112122\underline 112112122122112\cdots\\
\indent \Delta^1(w)= 21221121122122121121122121121221\cdots\\
\indent \Delta^2(w)=112212212112122112112\cdots\\
\indent \Delta^3(w)=2212112112212\cdots\\
\indent \Delta^4(w)=21121221\cdots\\
\indent \Delta^5(w)=12112\cdots\\
\indent \Delta^6(w)=112\cdots\\
\indent \Delta^7(w)=2\cdots$\\
$w$ has the factor $1121121$.

\item [(ii)] If $b \neq 2$, then\\
\indent $\displaystyle \Delta^0(w)= ((1^bb)^{2n+1}((1^bb^b)^{2n+1}(1b)^{2n+1})^{2n+1})^{2n+1}1^b((b1)^{2n+1}(b^b1^b)^{2n+1})^{2n+1}\\
\indent \qquad \qquad(b((1^bb^b)^{2n+1}(1b)^{2n+1})^{2n+1}1^b((b1)^{2n+1}(b^b1^b)^{2n+1})^{2n+1})^n(b(1^bb^b)^{2n+1}(1(b^b1^b)^{2n+1}b\\
\indent \qquad \qquad(1^bb^b)^{2n+1})^n(1b^b)^{2n+1}(1b)^{2n+1}((1^bb)^{2n+1}(1^bb^b)^{2n+1}(1b^b)^{2n+1}(1b)^{2n+1})^n((1^bb)^{2n+1}\\
\indent \qquad \qquad((1^bb^b)^{2n+1}(1b)^{2n+1})^{2n+1})^{2n+1} 1^b((b1)^{2n+1}(b^b1^b)^{2n+1})^{2n+1}(b((1^bb^b)^{2n+1}(1b)^{2n+1})^{2n+1 }1^b\\
\indent \qquad \qquad ((b1)^{2n+1}(b^b1^b)^{2n+1})^{2n+1})^n)^nb(1^bb^b)^{2n+1}(1(b^b1^b)^{2n+1}b(1^bb^b)^{2n+1})^n(1b^b)^{2n+1}((1b)^{2n+1}1^b\\
\indent \qquad \qquad ((b1)^{2n+1}b^b(1b)^{2n+1}1^b)^n((b1)^{2n+1}(b^b1^b)^{2n+1})^{2n+1}
b(1^bb^b)^{2n+1}(1(b^b1^b)^{2n+1}\\
\indent \qquad \qquad b(1^bb^b)^{2n+1})^n(1b^b)^{2n+1})^n(1b)^{2n+1}1^b((b1)^{2n+1}b^b(1b)^{2n+1}1^b)^n(((b1)^{2n+1}(b^b\underline 1^b)^{2n+1})^{2n+1}\\
\indent \qquad \qquad (b1^b)^{2n+1})^{2n+1} (b1)^{2n+1}b^b((1b)^{2n+1}1^b(b1)^{2n+1}b^b)^n(((1b)^{2n+1}(1^bb^b)^{2n+1})^{2n+1}(1b^b)^{2n+1})^{2n+1}\cdots\\
\indent \Delta^1(w)=((b1)^{2n+1}(b^b1^b)^{2n+1})^{2n+1}b(1^bb^b)^{2n+1}(1(b^b1^b)^{2n+1}b(1^bb^b)^{2n+1})^n((1b^b)^{2n+1}(1b)^{2n+1}1^b\\
\indent \qquad \qquad ((b1)^{2n+1}b^b(1b)^{2n+1}1^b)^n((b1)^{2n+1}(b^b1^b)^{2n+1})^{2n+1}b(1^bb^b)^{2n+1}(1(b^b1^b)^{2n+1}b(1^bb^b)^{2n+1})^n)^n\\
\indent \qquad \qquad(1b^b)^{2n+1}(1b)^{2n+1}((1^bb)^{2n+1}(1^bb^b)^{2n+1}(1b^b)^{2n+1}(1b)^{2n+1}((1^bb)^{2n+1}((1^bb^b)^{2n+1}\\
\indent \qquad \qquad (1b)^{2n+1})^{2n+1})^{2n+1}1^b(b1)^{2n+1}b^b1\\
\indent \Delta^2(w)=(1^bb^b)^{2n+1}(1b^b)^{2n+1}((1b)^{2n+1}(1^bb)^{2n+1}(1^bb^b)^{2n+1}(1b^b)^{2n+1})^n(1b)^{2n+1}1^b((b1)^{2n+1}\\
\indent \qquad \qquad b^b(1b)^{2n+1}1^b)^n ((b1)^{2n+1}(b^b1^b)^{2n+1})^{2n+1}b1^bb\cdots\\
\indent \Delta^3(w)=b^b(1b)^{2n+1}(1^b(b1)^{2n+1}b^b(1b)^{2n+1})^n(1^bb)^{2n+1}(1^bb^b)^{2n+1}1b\cdots\\
\indent \Delta^4(w)=(b1^b)^{2n+1}(b1)^{2n+1}b^b1\cdots\\
\indent \Delta^5(w)=(1b)^{2n+1}1^bb\cdots\\
\indent \Delta^6(w)=1^bb\cdots\\
\indent \Delta^7(w)=b\cdots$\\
$w$ has the factor $1^b(b1^b)^{2n+1}(b1)^{2n+1}$.
\end{itemize}

\item [Case (27)] If $p=1b1bbb$, there are $2$ subcases to consider.
\begin{itemize}
\item [(i)] If $b\neq 2$, then\\
\indent $\displaystyle \Delta^0(w)=((1^bb)^{2n+1}((1^bb^b)^{2n+1}(1b)^{2n+1})^{2n+1})^{2n+1}1^b(b1)^{2n+1}(b^b(1b)^{2n+1}1^b(b1)^{2n+1})^n\\
\indent \qquad (((b^b1)^{2n+1}((b^b1^b)^{2n+1}(b1)^{2n+1})^{2n+1})^{2n+1}b^b(1b)^{2n+1}(1^b(b1)^{2n+1}b^b(1b)^{2n+1})^n\\
\indent \qquad \qquad ((\underline 1^bb)^{2n+1}((1^bb^b)^{2n+1}(1b)^{2n+1})^{2n+1})^{2n+1}1^b(b1)^{2n+1}(b^b(1b)^{2n+1}1^b(b1)^{2n+1})^n)^n\\
\indent \qquad \qquad(b^b1)^{2n+1}(b^b1^b)^{2n+1}b1^bb\cdots\\
\indent \Delta^1(w)=(((b1)^{2n+1}(b^b1^b)^{2n+1})^{2n+1}(b1^b)^{2n+1})^{2n+1}(b1)^{2n+1}b^b1b\cdots\\
\indent \Delta^2(w)=((1^bb^b)^{2n+1}(1b)^{2n+1})^{2n+1}1^bb1\cdots\\
\indent \Delta^3(w)=(b^b1^b)^{2n+1}b1\cdots\\
\indent \Delta^4(w)=b^b1\cdots\\
\indent \Delta^5(w)=b\cdots$\\
$w$ has the factor $(( 1^bb)^{2n+1}((1^bb^b)^{2n+1}(1b)^{2n+1})^{2n+1})^{2n+1}1^b(b1)^{2n+1}$\\
\indent \qquad \qquad \qquad \qquad \qquad $(b^b(1b)^{2n+1}1^b(b1)^{2n+1})^n)^n(b^b1)^{2n+1}(b^b1^b)^{2n+1}b1^b$. 

\item [(ii)] If $b=2$, there are again 3 subcases to consider.
\begin{itemize}
\item [(a)] If $p=12122211$, then\\
\indent $\displaystyle \Delta^0(w)=112112212112122122\underline 11211221211212211\cdots\\
\indent \Delta^1(w)=21221121121221221121122\cdots\\
\indent \Delta^2(w)=112212112122122\cdots\\
\indent \Delta^3(w)=2211211212\cdots\\
\indent \Delta^4(w)=221211\cdots\\
\indent \Delta^5(w)=2112\cdots\\
\indent \Delta^6(w)=12\cdots\\
\indent \Delta^7(w)=1 \cdots$\\
$w$ has the factor $11211221211212211$.

\item [(b)] If $p=12122212$, then\\
\indent $\displaystyle \Delta^0(w)= 112112212112122122\underline 1121122121121122122\cdots\\
\indent \Delta^1(w)=212211211212212211212212\cdots\\
\indent \Delta^2(w)=11221211212211211\cdots\\
\indent \Delta^3(w)=22112112212\cdots\\
\indent \Delta^4(w)=2212211\cdots\\
\indent \Delta^5(w)=2122\cdots\\
\indent \Delta^6(w)=112\cdots\\
\indent \Delta^7(w)=2 \cdots$\\
$w$ has the factor $11211221211211$.

\item [(c)] If $p=1212222$, then\\
\indent $\displaystyle \Delta^0(w)=112112212112122122\underline 11211212211\cdots\\
\indent \Delta^1(w)=2122112112122121122\cdots\\
\indent \Delta^2(w)=1122121121122\cdots\\
\indent \Delta^3(w)=22112122\cdots\\
\indent \Delta^4(w)=22112\cdots\\
\indent \Delta^5(w)=221\cdots\\
\indent \Delta^6(w)=2 \cdots$\\
$w$ has the factor $1121121$. 
\end{itemize}

\end{itemize}

\item [Case (28)] If $p=1bb11b$, then\\
\indent $\displaystyle \Delta^0(w)=(1^bb^b)^{2n+1}(1b^b)^{2n+1}((1b)^{2n+1}(1^bb)^{2n+1}(1^bb^b)^{2n+1}(1b^b)^{2n+1})^n(1b)^{2n+1}1^b\\
\indent \qquad \qquad((b1)^{2n+1}b^b(1b)^{2n+1}1^b)^n((b1)^{2n+1}(b^b\underline 1^b)^{2n+1})^{2n+1}b1^bb\cdots\\
\indent \Delta^1(w)=b^b(1b)^{2n+1}(1^b(b1)^{2n+1}b^b(1b)^{2n+1})^n(1^bb)^{2n+1}(1^bb^b)^{2n+1}1b\cdots\\
\indent \Delta^2(w)=(b1^b)^{2n+1}(b1)^{2n+1}b^b1\cdots\\
\indent \Delta^3(w)=(1b)^{2n+1}1^bb\cdots\\
\indent \Delta^4(w)=1^bb\cdots\\
\indent \Delta^5(w)=b\cdots$\\
$w$ has the factor $1^bb1^b$.

\item [Case (29)] If $p=1bb1b$, then\\
\indent $\displaystyle \Delta^0(w)=(1^bb^b)^{2n+1}1((b^b1^b)^{2n+1}b(1^bb^b)^{2n+1}1)^n((b^b\underline 1^b)^{2n+1}(b1)^{2n+1})^{2n+1}b^b1b\cdots\\
\indent \Delta^1(w)=(b^b1)^{2n+1}(b^b1^b)^{2n+1}b1\cdots\\
\indent \Delta^2(w)=(b1)^{2n+1}b^b1\cdots\\
\indent \Delta^3(w)=1^bb\cdots\\
\indent \Delta^4(w)=b\cdots$\\
$w$ has the factor $1^bb1$.
\end{itemize}

\end{document}